# Adaptive and various learning-based algorithm for supply chain network equilibrium problems

Sheng-xue He

(Business School, University of Shanghai for Science and Technology, Shanghai 200093, China)

**Abstract**: This paper proposes a novel nonlinear programming model to capture the equilibrium state of complex supply chain networks. The model, equivalent to a variational inequality model, relaxes traditional strict assumptions to accommodate real-world complexities, such as nonlinear, non-convex, and non-smooth relationships between production, consumption, and pricing. To efficiently solve this challenging problem, we introduce a novel heuristic algorithm, the adaptive and various learning-based algorithm (AVLA), inspired by group learning behaviors. AVLA simulates individual learning processes within different subgroups at various stages, employing a success history-based parameter adaptation mechanism to reduce manual tuning. Extensive computational experiments on 29 benchmark problems and 5 supply chain networks demonstrate AVLA's superior performance compared to 19 state-of-the-art algorithms. AVLA consistently achieves the best results in terms of both average and best objective function values, making it a powerful tool for addressing complex supply chain network equilibrium problems.

## 1. Introduction

Supply chain networks (SCNs) comprise multiple decision-makers who aim to maximize their individual benefits through strategic interactions. This interplay often leads to an equilibrium state where no player can unilaterally improve their position. The resulting equilibrium quantities and prices guide decisions on production levels and pricing strategies, reflecting the relative power of each player within the network. A well-functioning SCN operates in a stable manner, benefiting all participants under fair conditions. Despite significant research efforts, modeling and solving SCNs remains a challenging task. Existing models, often based on variational inequalities, rely on simplifying assumptions about production processes, pricing mechanisms, and network structures. While these models provide valuable insights, they may yield unrealistic solutions when applied to real-world scenarios. To address these limitations, we propose a new supply chain network equilibrium model (SCNEM) that incorporates abstract functions to represent complex production and pricing dynamics. This novel approach, combined with the adaptive and various learning-based algorithm (AVLA), enables us to tackle highly nonlinear, non-convex, and non-smooth optimization problems. Extensive computational experiments demonstrate AVLA's superior performance in solving SCNs, often providing feasible solutions where other methods fail.

This paper offers two primary contributions. First, we introduce a generalized nonlinear programming model to capture supply chain network equilibrium under complex constraints. This model's flexibility allows it to accommodate various supply chain configurations by substituting abstract functions with concrete ones that reflect real-world conditions. Second, we present a novel population-based heuristic algorithm, AVLA, to efficiently solve the proposed model. Extensive computational experiments demonstrate AVLA's superior performance compared to 19 state-of-the-art algorithms.

The remainder of the paper is organized as follows. Section 2 reviews existing research on supply chain network equilibrium modeling and heuristic algorithms. Section 3 introduces the notation and formulates the generalized supply chain network equilibrium model, establishing its equivalence to a corresponding variational inequality model. Section 4 details the adaptive and various learning-based algorithm (AVLA), outlining its individual learning behaviors and parameter adaptation mechanism. Section 5 presents computational experiments on 29 benchmark problems and 5 supply chain network equilibrium problems, highlighting AVLA's effectiveness. Finally, Section 6 concludes the paper and suggests future research directions.

## 2. Literature review

### 2.1. Supply chain network equilibrium

Variational inequalities (VI) have been extensively used to model supply chain network equilibrium states. Common methods for solving the resulting variational inequality models (VIMs) include the Projection Method (PM), Modified Projection Method (MPM), and Euler Algorithm (EA). Nagurney et al. (2002) (Nagurney et al., 2002) were among the first to formulate a supply chain network with electronic commerce using the VIM framework. Huggins and Olsen (2003) (Huggins and Olsen, 2003) explored optimal policies for a two-stage supply chain under centralized control, where the downstream facility faced discrete stochastic demand. Subsequently, more complex supply chain network topologies were investigated by Dong et al. (2004) (Dong et al., 2004) and Dong et al. (2005) (Dong et al., 2005), incorporating random market demands and diverse decision-making objectives such as profit maximization, production efficiency, service level improvement, and transportation time minimization. In addition to these objectives, risk minimization has also been recognized as a critical criterion in multitiered supply chain network equilibrium models (Nagurney et al., 2005). Given the capability of VIs to capture interactions among decision-makers, this modeling approach has been extended to various network equilibrium problems, including multitiered financial networks (Nagurney and Ke, 2006), electric power supply chain networks (Wu et al., 2006) (Nagurney et al., 2007), and supply chains involving offshore outsourcing (Liu and Nagurney, 2011). Zhang (2006) (Zhang, 2006) introduced a general VIM to address competition among heterogeneous supply chains. Furthermore, the closed-loop supply chain network and coordination among manufacturers in multitiered supply chains were analyzed by Hammond and Beullens (2007) (Hammond and Beullens, 2007) and Hsueh and Chang (2008) (Hsueh and Chang, 2008), respectively. To solve their model, Hsueh and Chang (2008) (Hsueh and Chang, 2008) developed a method combining the diagonalization technique with user-equilibrium traffic assignment principles.

Over the past decade, research has expanded on earlier studies to address various specific fields within the supply chain. Methods have been adapted to tackle global outsourcing concerns and to create analytical frameworks that examine the impact of financial risks, utilizing specific mathematical approaches for model solutions(Liu and Nagurney, 2013) (Liu and Cruz, 2012). An equilibrium model has been introduced to assess sustainable fashion supply chains, while modeling techniques have been refined to address food supply chains, particularly considering the perishability of fresh produce(Nagurney and Yu, 2012) (Yu and Nagurney, 2013). Additionally, a nonlinear complementarity equilibrium model has been developed for supply chain networks using smoothing algorithms for resolution(Zhang and Zhou, 2012). An evolutionary equilibrium model has been created for closed-loop supply chains that account for seasonal demand variations(Feng et al., 2014). New frameworks have emerged that analyze supply chain networks over multiple periods while incorporating elements of opportunism risk(Daultani et al., 2015). The blood banking system has also been explored through the lens of supply chain networks, featuring models that investigate the implications of mergers and acquisitions within this domain(Masoumi et al., 2017). Stochastic robust optimization has been applied to closed-loop supply chain networks, highlighting the complexities involved in efficient supply chain management(Jabbarzadeh et al., 2018). Competitive models have been designed for food supply chains, incorporating iterative processes to refine outcomes(Nagurney et al., 2018). Recent developments include multi-period equilibrium models that adjust parameters dynamically based on loss-aversion, as well as extensions that consider behaviors across multiple attributes(Zhou et al., 2018) (Chan et al., 2019). Robust supply chain network models have been formulated with strategic decision-making principles, and newer approaches have integrated risk hedging via future contracts(Hirano and Narushima, 2019) (Liu and Wang, 2019). Lastly, dual-channel network models have been proposed for complex supply chain interactions, and considerations around labor shortages have been integrated into structural supply chain models)(Zhang et al., 2020) (Nagurney, 2021).

In recent years, sustainability, randomness, global outsourcing, and competition between supply chains have become significant areas of focus for researchers. An equilibrium model has been developed to analyze automotive supply chains, particularly considering uncertainties related to payment delays(Chen et al., 2020). Additionally, competition models have been established to examine how decisions made by suppliers and retailers affect wholesale pricing(Korpeoglu et al., 2020). Models integrating discrete choice and equilibrium modeling theory have also been introduced to explore time-cost competition between supply chains(Ma et al., 2020). An international human migration network model has been formulated as a variable income problem, addressing regulatory impacts on migration(Nagurney and Daniele, 2021). Furthermore, strategies for managing multiple supply chains toward sustainability have been explored(Brandao and Godinho-Filho, 2022), along with research on supply chain resilience in response to various disruptions(Carvalho et al., 2022), leading to the proposal of a resilience assessment index. The influence of trade policies on global supply chains has been investigated, emphasizing competitive dynamics within supply chains(Feng et al., 2022). Recent approaches include game-theoretic models that take into account uncertain demand and investments in green technology(Gupta et al., 2023), as well as analyses of dual-channel supply chains, focusing on optimal pricing strategies for manufacturers and retailers over multiple periods(He et al., 2023). Stochastic mixed integer programming models have been proposed to explore reshoring practices(Sawik, 2023), while complex network theory has been applied to understand the evolution of supply chains under random and targeted disruptions(Wang et al., 2023). Additionally, the introduction of buyback contracts in supply chain network modeling aims to balance after-sale services with product quality, providing a comprehensive framework for enhancing supply chain effectiveness(Xiao and Zhang, 2023). A four-tiered supply chain network was formulated into VIM by He (2023)(He and Cui, 2023c) where the strict limits to related functions are removed.

2.2. Heuristic algorithms

In the heuristic optimization field, researchers have proposed a huge number of algorithms in the past. Giving every algorithm a brief introduction is impossible in a paper due to the limited space. Each one of the existing heuristic algorithms usually provides a feasible way of solving all kinds of optimization problems and performs well on some special type of problems. Heuristic algorithms can be grouped into four categories roughly. In the following, we will introduce only part of them especially the original metaheuristic algorithms not the various later-developed variants.

The first category includes all kinds of evolution algorithms. Many advanced and powerful and widely used heuristic algorithms come from this category. Genetic Algorithm (GA) (Holland, 1992) may be the first widely accepted evolution algorithm. Later, algorithms such as Differential Evolution (DE) (Storn and Price, 1997) and the Covariance Matrix Adaptation Evolution Strategy(CMA-ES) (Hansen and Ostermeier, 1996) were introduced and successfully applied in practice. Tons of refined algorithms based on the earlier versions such as the CMA-ES including restarts with increasing population size (Beyer and Sendhoff, 2008) and the Combat Genetic Algorithm (CGA) (Auger and Hansen, 2005; Erol and Eksin, 2006) were developed. Heuristic algorithms developed from the DE obtained huge success in practice due to their simple implementation processes and the capacity of obtaining satisfying optimization results. Algorithms such as JADE (Jingqiao and Sanderson, 2007), Biogeography-Based Optimization (BBO) (Simon, 2008), the Success-History based parameter Adaptation for Differential Evolution (SHADE) (Tanabe and Fukunaga, 2013), L-SHADE (Tanabe and Fukunaga, 2014), the Backtracking Search Optimization Algorithm (BSO) (Civicioglu, 2013) and Symbiotic Organism Search (SOS) (Cheng and Prayogo, 2014) are just a few of them well-known to researchers.

The second category includes all kinds of swarm intelligent algorithms. One widely used swarm heuristic algorithm is the Ant Colony Optimization (ACO) presented by Dorigo in 1991 (Dorigo et al., 2006). ACO is especially easy to be applied to discrete optimization problems. Dolphin Echolocation Optimization (DEO) proposed by Kaveh and Farhoudi in 2013 (Kaveh and Farhoudi, 2013) is also fit for discrete optimization. Other two earlier developed and widely used swarm intelligent algorithms are Particle Swarm Optimization (PSO) (Kennedy and Eberhart, 1995) and Artificial Bee Colony algorithm (ABC) (Karaboga, 2005). Animals and plants provide rich sources for researcher to design metaheuristic algorithms. Many of the algorithms achieve surprising success in a wide area. The following algorithms are inspired by the behaviors of all kinds of animals: Shuffled Frog-Leaping Algorithm (SFLA) (Eusuff et al., 2006), Honey-Bee Mating Optimization (HBMO) (Afshar et al., 2007), Cuckoo Search with Levy Flights (Yang and Deb, 2009), Firefly Algorithm (Yang, 2010a), Bat Algorithm (Yang, 2010b) (Gandomi et al., 2013), Grey Wolf Optimizer (GWO) (Mirjalili et al., 2014), Coral Reefs Optimization Algorithm (CRO) (Salcedo-Sanz et al., 2014), Moth-flame optimization algorithm (MFO) (Mirjalili, 2015b), Ant Lion Optimizer (ALO) (Mirjalili, 2015a), Monkey King Evolution (MKE) (Meng and Pan, 2016), Whale Optimization Algorithm (WOA) (Mirjalili and Lewis, 2016), Dragonfly Algorithm (Mirjalili, 2016a), Grasshopper Optimization Algorithm (GOA) (Mirjalili et al., 2018), Salp Swarm Algorithm (SSA) (Mirjalili et al., 2017), Snake optimizer (SO) (Hashim and Hussien, 2022) and Giant Trevally Optimizer (GTO) (Sadeeq and Abdulazeez, 2022). In contrast to the algorithms inspired by animals' behaviors, algorithms such as Invasive Weed Optimization (IWO) (Mehrabian and Lucas, 2006), Plant Propagation Algorithm (PPA) (Salhi and Fraga, 2011) and Flower Pollination Algorithm (FPA) (Yang, 2012) are inspired by plant behaviors.

In addition to searching inspiring ideas from animals and plants, the physic and chemical phenomena and laws even math operations become the other rich area where many algorithms are inspired. Simulated annealing (Kirkpatrick et al., 1983) may be the first very successful algorithm come from this category. The following are some representatives including Big Bung-Big Crunch algorithm (BBBC) (Erol and Eksin, 2006), Central Force Optimization (CFO) (Formato, 2007), Gravity Search Algorithm (GSA) (Rashedi et al., 2009), Charged System Search (CSS) (Kaveh and Talatahari, 2010), Ray Optimization (Kaveh and Khayatazad, 2012), Water Cycle Algorithm (WCA) (Eskandar et al., 2012), Black Hole Algorithm (BHA) (Hatamlou, 2013), Global Mine Blast Algorithm (GMBA) (Sadollah et al., 2013; Yadav et al., 2020) (Yadav et al., 2020), Radial Movement Optimization (RMO) (Rahmani and Yusof, 2014), Enhanced Colliding Bodies Optimization (ECBO) (Kaveh and Ghazaan, 2014), Stochastic Fractal Search (SFS) (Salimi, 2015), Sine Cosine Algorithm (SCA) (Mirjalili, 2016b), Water Evaporation Optimization (WEO) (Kaveh and Bakhshpoori, 2016), Multi-Verse Optimizer (MVO) (Mirjalili et al., 2016), Vibrating Particles System algorithm (VPS) (Kaveh and Ilchi Ghazaan, 2017), Tug of War Optimization (TWO) (Kaveh and Zolghadr, 2016), Thermal Exchange Optimization (TEO) (Kaveh and Dadras, 2017), Opposition-based tug of war optimization 2020 (Nguyen et al., 2020), and Equilibrium Optimizer (EO) (Faramarzi et al., 2020).

The last category includes metaheuristic algorithms inspired by phenomena or actions associated with human society. Generally, Tabu Search (Glover, 1990) is viewed as the first one in this category. Harmony Search (HS) (Lee and Geem, 2005) is another algorithm in this category widely used in practice. Other algorithms such as Teaching-Learning-Based Optimization (TLBO) (Rao et al., 2011), Gaining-Sharing Knowledge based algorithm (GSK) (Mohamed et al., 2020), Life Choice-Based Optimizer (LCBO) (2019) (Khatri et al., 2020) and Medalist Learning Algorithm (MLA) (He and Cui, 2023a) (He and Cui, 2023b) (He, 2023)base their main idea on the learning activities of human being. Other representatives include Imperialist Competitive Algorithm (ICA) (Atashpaz-Gargari and Lucas, 2007), League Championship Algorithm (LCA) (Kashan, 2009), Jaya (Rao, 2019) and Levy flight Jaya algorithm (LJA) (Iacca et al., 2021). Among them, the implementation process of Jaya is very simple. A recently proposed algorithm in this category is Hiking Optimization Algorithm (HOA) (Oladejo et al., 2024). The new algorithm-AVLA to be presented in this paper is also included in this category. AVLA is different from the other existing learning-based algorithm in this category by the concrete forms of formulating the learning behaviors.

## 3. Supply chain network equilibrium problem

### 3.1. Main notations

In a supply chain network (SCN), we will call any decision maker a spot (or node) usually denoted by $i$ and $j$ in the subsequent related expressions. Among the decisionmakers, the set of all the suppliers is denoted by $I_S$. The set of all markets is denoted by $I_M$. A market is also viewed as a spot in SCN. Except suppliers and markets, the others include manufacturer, wholesaler, and retailer are put in a common set $I_{mwr}$. The set of all the decision makers is defined as $I \equiv I_S \cup I_{mwr} \cup I_M$.

If there exists any transaction between two decision makers, we say there is a link connecting them. A link denoted by $(i,j)$ is the transaction connector from spot $i$ to spot $j$. The set of all the products including raw material, component of product, half made product, and finished product is denoted by $PK$. A typical product is denoted by $k \in PK$. The quantity of product $k$ on link $(i,j)$ is denoted by $f_{i,j,k}$. The upper bound to product $k$ on link $(i,j)$ is denoted by $f_{i,j,k}^{Max}$. We call the flow vector composed of all the link flows $f_{i,j,k}, \forall i,j,k$ a flow pattern which satisfies the flow conservation conditions at all the spots of a supply chain network.

At a spot $i$, we will distinguish the quantities of products as follows. $Q_{i,k}^E$, $Q_{i,k}^H$ and $Q_{i,k}^O$ are the arrived, holding and sold products $k$ at $i$, respectively. $Q_i^E$, $Q_i^H$ and $Q_i^O$ are the vectors of $Q_{i,k}^E$, $Q_{i,k}^H$ and $Q_{i,k}^O$, $k \in PK$ at $i$, respectively. $Q_i$ is the vector of products composed of $Q_i^E$, $Q_i^H$ and $Q_i^O$. $Q_{i,k}^{E,Max}$ denotes the maximal quantity of raw material of $k \in PK_S$ which could be supplied by a supplier $i \in I_S$. $PK_S$ is the set of raw material.

Notations about prices of products are defined as follows. $p_{i,k}^O$ and $p_{i,k}^E$ are the selling and buying prices of product $k \in PK$ at a spot $i$, respectively. For a product $k \in PK$ which comes from spot $j$ by link $(j,i)$, $p_{j,i,k}^E$ denotes the possible buying price of product $k \in PK$ at spot $i$ offered by decision maker $j$.

The cost associated with shipping $f_{j,i,k}$ amount of product $k$ from $j$ to $i$ is denoted by $C_{j,i,k}(f_{j,i,k})$. The total cost incurred by product $k$ at $i$ is denoted by $C_{i,k}(Q_i)$ which generally includes the related purchase expenditure, manufactural cost, holding cost, and transaction cost.

We assume that all decision makers expect a positive profit rate. Let $\lambda_{i,k}$ be the profit rate of product $k$ for decision maker $i$ and $\lambda_{i,k}^{Max}$ the corresponding maximal profit rate.

## 3.2 Supply chain network equilibrium model (SCNEM)

With the notations given in the preceding subsection, we define $PA_{i,j,k} = max\{0, p_{i,k}^O + C_{i,j,k}(f_{i,j,k}) - p_{j,k}^E\}$ and $PB_{i,j,k} = max\{0, -p_{i,k}^O - C_{i,j,k}(f_{i,j,k}) + p_{j,k}^E\}$. Now we can formulate the supply chain network equilibrium model (SCNEM) as follows:

$$\min \sum_{i,j,k} \{f_{i,j,k} PA_{i,j,k} + (f_{i,j,k}^{Max} - f_{i,j,k}) PB_{i,j,k}\} \tag{1}$$

$$\sum_i f_{i,j,k} = Q_{j,k}^E, \quad \forall j,k \tag{2}$$

$$Q_{j,k}^O = \sum_i f_{j,i,k}, \quad \forall j,k \tag{3}$$

$$[Q_i^O, Q_i^H] = M(Q_i^E), \quad \forall i \in I \backslash I_M \tag{4}$$

$$p_{i,k}^O = \frac{C_{i,k}(Q_i)}{Q_{i,k}^O + Q_{i,k}^H}(1 + \lambda_{i,k}), \forall i \in I \backslash I_M, k \tag{5}$$

$$p_{j,i,k}^E = p_{j,k}^O + C_{j,i,k}(f_{j,i,k}), \quad \forall j,i,k \tag{6}$$

$$p_{i,k}^E = \min_j \{p_{j,i,k}^E\}, \forall i \in I \backslash I_M, k \tag{7}$$

$$p_{i,k}^E = P(Q_{i,k}^E), \forall i \in I_M \tag{8}$$

$$0 \leq f_{i,j,k} \leq f_{i,j,k}^{Max}, \forall i,j,k \tag{9}$$

$$0 \leq \lambda_{i,k} \leq \lambda_{i,k}^{Max}, \forall i \in I \backslash I_M, k \tag{10}$$

$$0 \leq Q_{i,k}^E \leq Q_{i,k}^{E,Max}, \quad \forall i \in I_S, k \tag{11}$$

The objective function given in equation (1) comes from the following inequalities:

$$(p_{j,k}^{O*} + C_{j,i,k}(f_{j,i,k}^*) - p_{i,k}^{E*}) \times (f_{j,i,k} - f_{j,i,k}^*) \geq 0, \forall i,j,k \tag{12}$$

The inequality (12) states that if the offered price of product $k$ at spot $i$ offered by spot $j$ is greater than the actual buying price of product $k$ at spot $i$, there is no actual transaction between $i$ and $j$. Conversely, if the offered price of product $k$ at spot $i$ offered by spot $j$ is lower than the actual buying price of product $k$ at spot $i$, the actual transaction of product $k$ between $i$ and $j$ will occur at the highest level $f_{i,j,k}^{Max}$. Here $p_{j,k}^{O*} + C_{j,i,k}(f_{j,i,k}^*)$ is the price of product $k$ at spot $i$ offered by spot $j$ and $p_{i,k}^{E*}$ is the actual buying price of product $k$ at spot $i$. Theorem 1 will formally prove the equivalence between the nonlinear optimization problem (NP) and the variational inequality problem (VIP). The NP, defined by objective function (1) and constraints (2-11), is equivalent to the VIP, defined by inequality (12) and feasible space constraints (2-12).

Equations (2) and (3) are the flow conversation constraints of products. Equation (4) is an abstract manufacturing function which transforms the arrived products (including raw materials, parts and semi-finished products) into the product to be held or sold at spot $i$. The abstract form of manufacturing function $[Q_i^O, Q_i^H] = M(Q_i^E)$ will be replaced by any real transformation function when applying the above model in real life.

Equations (5) to (8) are the constraints about pricing. Equation (5) points out that the selling price $p_{i,k}^O$ of a product $k$ at spot $i$ is determined by the associated total cost $C_{i,k}(Q_i)$ and the expected profit rate $\lambda_{i,k}$. Equation (6) defines the offered price $p_{j,i,k}^E$ of product $k$ at spot $i$ come from spot $j$. Equation (7) defines the actual purchasing price $p_{i,k}^E$ of product $k$ at spot $i$ which is not a market. Equation (8) defines the specified $p_{i,k}^E$ at a market $i$. The abstract function $P(Q_{i,k}^E)$ can be replaced by any proper concrete one in practice.

Constraints (9-11) give the limits to the feasible ranges of variables.

The SCNEM composed of equations (1-11) is very general and covers all kinds of supply chain structures, such as the popular multilevel supply chain network. This SCNEM is highly nonlinear and non-smooth. The classic optimization tools which require derivative information are unsuitable for solving it. In the next section, we will present a powerful population-based heuristic algorithm for researchers and practitioners to solve this SCNEM.

**Theorem 1:** If the variational inequality problem (VIP) (2-12) is satisfied by a flow pattern $f^*$, then the nonlinear programming problem (NP) (1-11) reaches its minimum at $f^*$. If a feasible flow pattern $f^*$ makes the objective value of NP (1-11) equal 0, then it will satisfy the VIP (2-12). A flow pattern is the vector composed of all the link flows associated with a given state of SCN where the constraints (2-4) are satisfied.

**Proof:** The proof of the first statement of this theorem is presented as follows. One fact to be used later is that the objective function value of NP (1) is greater than and equal to 0 because any link flow needs to satisfy the constraint (9). If the inequality (12) holds, three situations may arise.

In the first situation, if $p_{j,k}^{O*} + C_{j,i,k}(f_{j,i,k}^*) - p_{i,k}^{E*} > 0$, then $f_{j,i,k} - f_{j,i,k}^* \geq 0$ for all $0 \leq f_{j,i,k} \leq f_{j,i,k}^{Max}$. In other words, $p_{j,k}^{O*} + C_{j,i,k}(f_{j,i,k}^*) - p_{i,k}^{E*} > 0$ leads to $f_{j,i,k}^* = 0$. Obviously, in this case, the corresponding items in (1), that is $f_{j,i,k} PA_{j,i,k} + (f_{j,i,k}^{Max} - f_{j,i,k}) PB_{j,i,k}$, equals 0.0 when replacing the price and flow with the corresponding equilibrium values, such as replacing $f_{j,i,k}$ with $f_{j,i,k}^*$.

In the second situation, if $p_{j,k}^{O*} + C_{j,i,k}(f_{j,i,k}^*) - p_{i,k}^{E*} < 0$, then $f_{j,i,k} - f_{j,i,k}^* \leq 0$ for all $0 \leq f_{j,i,k} \leq f_{j,i,k}^{Max}$. In other words, $p_{j,k}^{O*} + C_{j,i,k}(f_{j,i,k}^*) - p_{i,k}^{E*} < 0$ leads to $f_{j,i,k}^* = f_{j,i,k}^{Max}$. Similarly, in this case, the corresponding items in (1), that is $f_{j,i,k} PA_{j,i,k} + (f_{j,i,k}^{Max} - f_{j,i,k}) PB_{j,i,k}$, equals 0.0 when replacing the price and flow with the corresponding equilibrium values, such as replacing $f_{j,i,k}$ with $f_{j,i,k}^*$.

In the third situation, if $p_{j,k}^{O*} + C_{j,i,k}(f_{j,i,k}^*) - p_{i,k}^{E*} = 0$, then any feasible $f_{j,i,k}^* \in [0, f_{j,i,k}^{Max}]$ makes items $f_{j,i,k} PA_{j,i,k} + (f_{j,i,k}^{Max} - f_{j,i,k}) PB_{j,i,k}$ equals 0 when replacing the price and flow with the corresponding equilibrium values, such as replacing $f_{j,i,k}$ with $f_{j,i,k}^*$.

In view of that the objective function value is greater than and equal to 0, the above three situations prove that if a flow pattern satisfies the VIP (2-12), the NP (1-11) must achieve its minimum value with this flow pattern.

The proof of the second statement of this theorem is presented as follows. Let us check if $f_{j,i,k}PA_{j,i,k} + (f_{j,i,k}^{Max} - f_{j,i,k})PB_{j,i,k}$ equals 0 for a given flow pattern $f_{j,i,k}^*$, $\forall j,i,k$, the corresponding inequality $(p_{j,k}^{O*} + C_{j,i,k}(f_{j,i,k}^*) - p_{i,k}^{E*}) \times (f_{j,i,k} - f_{j,i,k}^*) \geq 0, \forall i,j,k$ must hold. There are also three situations required to be considered here.

In the first situation, if $f_{j,i,k}^*=0$, then $(f_{j,i,k}^{Max} - f_{j,i,k})max\{0, -p_{j,k}^O - C_{j,i,k}(f_{j,i,k}) + p_{i,k}^E\}$ must equal 0 because $f_{j,i,k}max\{0, p_{j,k}^O + C_{j,i,k}(f_{j,i,k}) - p_{i,k}^E\}$ equals 0 and the whole expression is assumed to be 0. In this case, we can further deduce that $p_{i,k}^E - p_{j,k}^O - C_{j,i,k}(f_{j,i,k})$ must be less than or equal to 0. In other words, $(p_{j,k}^{O*} + C_{j,i,k}(f_{j,i,k}^*) - p_{i,k}^{E*}) \times (f_{j,i,k} - f_{j,i,k}^*) \geq 0$ because $f_{j,i,k} - f_{j,i,k}^* = f_{j,i,k} \geq 0$ and $p_{j,k}^{O*} + C_{j,i,k}(f_{j,i,k}^*) - p_{i,k}^{E*} \geq 0$.

In the second situation, if $f_{j,i,k}^* = f_{j,i,k}^{Max}$, then $f_{j,i,k}max\{0, p_{j,k}^O + C_{j,i,k}(f_{j,i,k}) - p_{i,k}^E\}$ must equal 0 because $(f_{j,i,k}^{Max} - f_{j,i,k})max\{0, -p_{j,k}^O - C_{j,i,k}(f_{j,i,k}) + p_{i,k}^E\}$ equals 0 and the whole expression is assumed to be 0. In this case, we can further deduce that $p_{j,k}^O + C_{j,i,k}(f_{j,i,k}) - p_{i,k}^E$ must be less than or equal to 0. In other words, $(p_{j,k}^{O*} + C_{j,i,k}(f_{j,i,k}^*) - p_{i,k}^{E*}) \times (f_{j,i,k} - f_{j,i,k}^*) \geq 0$ because $f_{j,i,k} - f_{j,i,k}^* = f_{j,i,k} - f_{j,i,k}^{Max} \leq 0$ and $p_{j,k}^O + C_{j,i,k}(f_{j,i,k}) - p_{i,k}^E \leq 0$.

In the third situation if $0 < f_{j,i,k}^* < f_{j,i,k}^{Max}$, then $p_{i,(i,j)} + c_{ij}(f_{ij}) - p_{j,(i,j)}$ must equal 0.0 because we assume that $f_{j,i,k}max\{0, p_{j,k}^O + C_{j,i,k}(f_{j,i,k}) - p_{i,k}^E\} + (f_{j,i,k}^{Max} - f_{j,i,k})max\{0, -p_{j,k}^O - C_{j,i,k}(f_{j,i,k}) + p_{i,k}^E\}$ equals 0. In this case, it is easy to see that $(p_{j,k}^{O*} + C_{j,i,k}(f_{j,i,k}^*) - p_{i,k}^{E*}) \times (f_{j,i,k} - f_{j,i,k}^*) \geq 0$.

By summarizing the above three situations, we can conclude that if a feasible flow pattern $f^*$ makes the objective value of NP (1-11) equal to 0, then it will make VIP (2-12) hold. □

## 4. Designing the adaptive and various learning-based algorithm

In this section, we will first introduce the notations used to describe a general optimization problem. Then we will formulate the learning behaviors of individuals into mathematical operations that are the main operations of AVLA. Subsequently, the way to adapt parameters in AVLA is introduced. At the end, the implementation process of AVLA is summarized.

### 4.1. General optimization problem description

Though our main goal is to use AVLA to solve SCNEM given in the preceding section, AVLA can be applied to general optimization problems. We can view SCNEM as a special case of general optimization problem. In the following, we can introduce the basic notations used to describe general optimization problems.

To solve an optimization problem, we use $X \coloneqq \{x_i | i = 1,2, ..., N\}$ to denote the solution set with $x_i$ as its representative. In AVLA, a solution is a member in the learning group. Here $N$ is the size of the population. A solution is composed of $D$ elements. Let $x_i = (x_{i,1}, x_{i,2}, ..., x_{i,D})$. We assume that any element needs to be chosen from a given set. That is $x_{i,j} \in V_j$ for all $i$ where $V_j$ is the feasible range or set of the $j$th element of any member. Assume that $\alpha_j \leq x_{i,j} \leq \beta_j$, $\forall i$. The objective function value (or the penalty function value corresponding to a constrained optimization problem) with respect to member $x_i$ is denoted by $fit(x_i)$. We aim to minimizing the objective function value. The initial solutions are generated randomly in the feasible space. For a member $x_i$, its initial $j$th element $x_{i,j}$ is generated as follows:

$$x_{i,j} = \alpha_j + rand[0,1](\beta_j - \alpha_j) \tag{13}$$

where $rand[0,1]$ is random number picked from interval $[0,1]$ according to the uniform distribution.

## 4.2. Learning process modeling

In this subsection we will model the learning actions of members in a learning group. Firstly, we will deal with the learning behaviors of members with high performance. We will call a member with high performance an elite in this paper. Secondly, the learning actions of common learners will be formulated into mathematic functions. At the end, we will model the reflection actions of members with unsatisfied performance and the whole group, respectively.

### *4.2.1. Learning among elites*

Confucius said: "Let me take a stroll with any two people, and I can always be sure of learning something from them. I can take their good points and emulate them, and I can take their bad points and correct them in myself." The elites in a learning group can put the above saying in practice as follows. At first, an elite can choose two other elites randomly as his/her models to emulate their strategies; then by comparing performances, the elite can determine whether to follow or to keep away from the selected elites.

In view of that the actual learning efficiency is always influenced by all kinds of factors and shows obvious uncertainty, we assume that the learning of an individual will be realized through three steps. Three steps are the ideal learning, practical learning, and actual learning. They are modeled in equations (14), (15) and (16), respectively, as given below:

$$v_{e,j}^t = x_{e,j}^t + s_1 F_e(x_{e1,j}^t - x_{e,j}^t) + s_2 F_e(x_{e2,j}^t - x_{e,j}^t) \tag{14}$$

$$v_{e,j}^{t+1} = \begin{cases} v_{e,j}^t & \text{if } rand[0,1] \leq CR_e \text{ or } j = j_{rand} \\ x_{e,j}^t & \text{otherwise} \end{cases} \tag{15}$$

$$x_e^{t+1} = \begin{cases} v_e^{t+1} & \text{if } fit(v_e^{t+1}) < fit(x_e^t) \\ x_e^t & \text{otherwise} \end{cases} \tag{16}$$

In above, $e$ indicates the elite in question; $e1$ and $e2$ are two randomly chosen different elites. $t$ is the learning stage that is corresponding to the iteration of AVLA. $s_1$ and $s_2$ are two sign parameters. If $fit(x_e^t) > fit(x_{e1}^t)$, $s_1$ equals 1; or else, $s_1$ equals -1. Similarly, if $fit(x_e^t) > fit(x_{e2}^t)$, $s_2$ equals 1; or else, $s_2$ equals -1. $F_e$ is *the learning acceptation rate* of member $e$. $F_e$ will be limited to [0, 1]. The rule of how to change $F_e$ will be introduced in **Section 4.3.1**. $v_{e,j}^t$ can be viewed as the ideal result from learning from other two elites because in this situation, $e$ can refine all his/her elements in a full scale.

In reality, it is nearly impossible for a member to learn everything from another member. A practical learning is to just emulate part of the features of the selected model. Here a feature stands for an element of a solution to the optimization problem in question. Equation (15) formulates the above thought with a stochastic judgement. $j_{rand}$ is an index chosen randomly from 1, 2, …, D. $rand[0,1]$ is a random number chosen from interval $[0,1]$ according to uniform distribution. $CR_e$ is *the practical learning factor*. The rule of how to change $CR_e$ will be introduced in **Section 4.3.1**.

We assume that an elite only accepts the change which can improve his/her current performance. Equation (16) models this decision-making behavior called the actual learning. The result from the above three steps will be the new state of elite $e$ to be estimated later.

### *4.2.2. Learning by commons*

Confucius said: "I listen a lot, pick the best of it, and follow it; I observe a lot and take note of it. This is the best way for me to learn." For a common learner whose performance is lower than any elite, the learning process is similar to the elites but with different learning models to emulate.

A common member stands for an individual other than elites in AVLA. We assume that a common member may choose to learn from an elite and one common member or decide to learn from two common members. The learning process followed by a common member consists of three steps similar to that related to the learning process among elites. The three steps composed of the ideal learning, practical learning and actual learning are formulated in equations (17), (18) and (19), respectively, as follows:

$$v_{i,j}^t = \begin{cases} x_{i,j}^t + s_1 F_i(x_{i1,j}^t - x_{i,j}^t) + s_2 F_i(x_{i2,j}^t - x_{i,j}^t) & if\ rand[0,1] > \text{LE}(t) \\ x_{i,j}^t + F_i(x_{e,j}^t - x_{i,j}^t) + s_2 F_i(x_{i2,j}^t - x_{i,j}^t) & otherwise \end{cases} \quad (17)$$

$$v_{i,j}^{t+1} = \begin{cases} v_{i,j}^t & if\ rand[0,1) \leq CR_i\ or\ j = j_{rand} \\ x_{i,j}^t & otherwise \end{cases} \quad (18)$$

$$x_i^{t+1} = \begin{cases} v_i^{t+1} & if\ fit(v_i^{t+1}) < fit(x_i^t) \\ x_i^t & otherwise \end{cases} \quad (19)$$

In equation (17), $i$ denotes the common learner in question; $e$ stands for a randomly selected elite; $i1$ and $i2$ are two members randomly chosen from the common members different from $i$ and $e$. $s_1$ and $s_2$ are two sign parameters. If $fit(x_i^t) > fit(x_{i1}^t)$, $s_1$ equals 1; or else, $s_1$ equals -1. Similarly, if $fit(x_i^t) > fit(x_{i2}^t)$, $s_2$ equals 1; or else, $s_2$ equals -1. $F_i$ and $CR_i$ are the learning acceptation rate and the practical learning factor of common member $i$. The explanation of equations (18)-(19) is similar to the explanation of equations (15)-(16).

In equation (17), LE(t) stands for the probability of that a common member decides to learn from an elite not another common member. LE(t) is defined by the following equation:

$$\text{LE}(t) = \left[1 + e^{\frac{2\gamma}{maxNumIter}(\frac{maxNumIter}{2} - t)}\right]^{-1} \quad (20)$$

In equation (20), $maxNumIter$ is the maximal number of iterations allowed to run AVLA once to solve an optimization problem. $t$ is the counter of the current iteration. $\gamma$ is called the truncated radius of learning curve. In this paper, we just assign 6 to $\gamma$. This curve function (20) is also used in MLA (He and Cui, 2023b). With the increase of $t$, the value of LE($t$) gradually increases from a low positive value near to zero to a high value a little bit less than 1. With LE(t) in (17), the common member is more likely to learn from other common learners at the beginning, but with the increase of iteration, he /she will be prone to learning from some elite.

### 4.2.3. Reflection after sorting

Confucius said: "Learning without thinking leads to perplexity. Thinking without learning leads to trouble." After the previous mentioned learning, members will be sorted from the best to the worst according to their current performances. We assume that members will conduct some forms of reflections to further refine their performance after the sorting.

Before presenting the concrete reflection way, we need the following notation. Let $x_{i,j}$ that is the $j$th element of the member $i$ be limited to the interval $[\alpha_j, \beta_j]$. We define $x_{i,j}^R$ the reflection position of $x_{i,j}$ in $[\alpha_j, \beta_j]$ such that:

$$x_{i,j}^R = \alpha_j + \beta_j - x_{i,j} \quad (21)$$

All the reflection positions of the elements of $x_i$ form a new spot in the searching space denoted by $x_i^R$. We call $x_i^R$ the *opposite position* of $x_i$ in the search space.

Now with the new conception of opposite position, we can formulate the reflection behavior of a member with unsatisfied performance as follows:

$$x_i = \begin{cases} x_i^R & \text{if } fit(x_i^R) < fit(x_i) \\ x_i^{STOCH} & \text{otherwise} \end{cases} \quad \forall\, N - n_E < i \leq N \tag{22}$$

In equation (22), $N$ and $n_E$ are the size of population and the size of elites, respectively. $x_i^{STOCH}$ stands for a solution generated randomly with equation (13). Equation (22) means that if a member's position is at the tail part of the list of sorted members, he/she will do reflection of trying his/her opposite position in the search space. If his/her opposite position has a better outcome, he/she will occupy the new position; or else, he/she will choose a random position, in other words being replaced by a new member.

In addition to the above reflection carrying out by members with unsatisfied performance, the whole group may conduct some reflection if the performance of the group as a whole is not improved for a long time. In AVLA, the performance of the group as a whole is indicated by the best-so-far objective function value. If the best-so-far objective function value has not been improved for a given number $n_R$ of successive iterations, all the members will conduct reflection as follows:

$$x_i = \begin{cases} x_i^R & \text{if } fit(x_i^R) < fit(x_i) \text{ and } i \leq N - n_E \\ x_i & \text{if } fit(x_i^R) \geq fit(x_i) \text{ and } i \leq N - n_E \\ x_i^R & \text{if } N - n_E < i \leq N \end{cases} \tag{23}$$

Equation (23) means that if a member who is not a member with unsatisfied performance, he/she will check his/her opposite position and occupy it if it can lead to better performance than his/her current position. But if the member has an unsatisfied performance, he/she will be forced to move to his/her opposite position.

### 4.3. Parameter adaptation

#### 4.3.1. Adapt the learning acceptation rate and the practical learning factor

To obtain the expected performance of a heuristic algorithm applied to a given problem, researchers generally need to adjust the parameters by trial-and-error. In general, there is no fixed parameter setting suitable for all problems and even the different phases during the process of optimization. To lighten the above burden and improve the resilience of a heuristic algorithm, parameter adaptation is adopted in many algorithms, e.g. JADE. The successful applications in the past in optimization field encourage the author to adopt it in AVLA. The adaptive method adopted in AVLA is the success-history based parameter adaptation first used in SHADE.

**Table 1** Historical memory of the learning acceptation rate and the practical learning factor.

| Index | 1 | 2 | ... | H-1 | H |
|---|---|---|---|---|---|
| $M_{CR}$ | $M_{CR,1}$ | $M_{CR,2}$ | ... | $M_{CR,H-1}$ | $M_{CR,H}$ |
| $M_F$ | $M_{F,1}$ | $M_{F,2}$ | ... | $M_{F,H-1}$ | $M_{F,H}$ |

AVLA maintains a historical memory with $H$ entries for both the learning acceptation rate $F$ and the practical learning factor $CR$ as shown in Table 1. The initial values of all the entries will be set to 0.5. At the beginning of each iteration of AVLA, every member $i$ in the learning group will select his/her learning acceptation rate $F_i$ and practical learning factor $CR_i$ as follows:

$$CR_i = randn_i(M_{CR,r_i}, 0.1) \tag{24}$$
$$F_i = randc_i(M_{F,r_i}, 0.1) \tag{25}$$

In equations (24) and (25), $r_i$ is a randomly selected index from $[1, H]$. $randn_i(M_{CR,r_i}, 0.1)$ means to generate $CR_i$ according to a normal distribution with mean $M_{CR,r_i}$ and standard deviation 0.1. If the generated $CR_i$ is outside of $[0,1]$, it will be truncated. $randc_i(M_{F,r_i}, 0.1)$ means to generate $F_i$ according to a Cauchy distribution with location parameter $M_{F,r_i}$ and scale parameter 0.1. If the generated $F_i$ is greater than 1, it will be replaced by 1; if the value is less than 0, the generation process using Cauchy distribution will be repeated until a proper value of $F_i$ is generated.

At the end of each iteration, the historical memory will be updated. During the learning process in an iteration, if a member $i$ successfully improves his/her performance, his/her learning acceptation rate $F_i$ and his/her practical learning factor $CR_i$ will be stored in sets $S_{CR}$ and $S_F$, respectively. The historical memory is updated as follows:

$$M_{CR,k,t+1} = \begin{cases} mean_{WL}(S_{CR}) & if\ S_{CR} \neq \emptyset \\ M_{CR,k,t} & otherwise \end{cases} \quad (26)$$

$$M_{F,k,t+1} = \begin{cases} mean_{WL}(S_F) & if\ S_F \neq \emptyset \\ M_{F,k,t} & otherwise \end{cases} \quad (27)$$

In equations (26) and (27), the subscript $k$ is the index of the memory entry to be updated. The initial value of $k$ is set to 1. Then whenever a new element is inserted into the memory, the value of $k$ will be increased by 1. If the value of $k$ becomes greater than $H$, it will be reset to 1.

The equations (26) and (27) use the weighted Lehmer means defined in equations (28) and (29):

$$mean_{WL}(S_{CR}) = \frac{\sum_{k=1}^{|S_{CR}|} w_k S_{CR,k}^2}{\sum_{k=1}^{|S_{CR}|} w_k S_{CR,k}} \quad (28)$$

$$mean_{WL}(S_F) = \frac{\sum_{k=1}^{|S_F|} w_k S_{F,k}^2}{\sum_{k=1}^{|S_F|} w_k S_{F,k}} \quad (29)$$

The weight $w_k$ is defined as follows:

$$w_k = \frac{\Delta fit_k}{\sum_{l=1}^{|S_{CR}|} \Delta fit_l} \quad (30)$$

$$\Delta fit_k = |fit(v_k^t) - fit(x_k^t)| \quad (31)$$

In equation (31), $\Delta fit_k$ is the absolute value of the change of the objective function values for member $k$ after the three-steps individual learning.

### 4.3.2. Increase the size of elites

The size of elites has a noticeable impact on the convergence of AVLA. At the earlier phase of search, the information held by elites is very likely uncertain and misleading regarding to the final best solution. Learning from them at the earlier phase of the search will do no obvious good to other members. In view of the above observation, only a few of members with relative high performance should be treated with as elites at the earlier phase of learning. When the optimization keeps going on, elites as a whole are more likely holding useful information of the direction of further refinement. At the same time, learning among elites can speed up the convergence to the promising position of the search space at the ending phase of learning. In view of the above observation, we will adapt the size of elites according to the following equation:

$$n_E(t) = round\left[3 + \frac{t(0.2N-3)}{maxNumIter}\right] \quad (32)$$

In equation (32), $t$ is the index of the current iteration and $maxNumIter$ is the maximal number of iterations predetermined for AVLA. It is easy to see that the initial size of elites is set to 3. With the increase of iterations, the size of elites will gradually be increased to $0.2N$ at the end.

The initial and final size of elites may be changed with respect to the problem in question. The above choice is just based on the experience of the author and it works most of time.

## 4.4. The implementation process of AVLA

In the following we sum up the implementation process of AVLA as given previously.

**Step 1 Initialization**

Initialize values for $N$, $maxNumIter$, $n_R$ and $H$. And set $nR = 0$.

**Step 2 Generate the initial learning group**

Generate the initial learning group $X \coloneqq \{x_i | i = 1,2,\dots,N\}$ using equation (13). And sort the members according to their performances.

**Step 3 Conduct learning**

   **Step 3.1** For every elite member, carry out the three-steps learning by using equations (14-16).

   **Step 3.2** For every common learner, carry out the three-steps learning by using equations (17-19).

   **Step 3.3** Sort the members according to their new performances.

**Step 4 Update historical memory**

Update the historical memory by equations (26) and (27).

**Step 5 Do necessary reflections**

If the best-so-far solution has not changed after $n_R$ successive iterations, in other words $nR$ equals $n_R$, set $nR = 0$ and execute **Step 5.1**; or else, increase $nR$ by 1 and carry out **Step 5.2**.

   **Step 5.1** Do reflection of the whole group using equation (23). Jump to **Step 5.3**.

   **Step 5.2** Do reflection by members with unsatisfied performance using equation (22).

   **Step 5.3** Sort the members according to their new relative performance.

**Step 6 Check terminal condition**

If the terminal criteria are satisfied, output the best-so-far solution and end the program; or else, go to **Step 3**.

## 5. Numerical experiments

### 5.1. Setting of numerical analysis

In this section, we apply the Adaptive and Various Learning-Based Algorithm (AVLA) and its non-adaptive variant, the Various Learning-Based Algorithm (VLA), to 29 well-known test benchmarks and 5 supply chain network equilibrium problems. These are compared against 18 other heuristic algorithms. Unlike AVLA, VLA omits parameter adaptation. Specifically, the learning acceptance rate in VLA is replaced by a random value drawn from a uniform distribution over [0,1], and the practical learning factor is replaced by a fixed constant.

Given the large number of existing heuristic algorithms, we selected representative ones with proven performance widely recognized in practice. The 18 algorithms compared include: Covariance Matrix Adaptation Evolution Strategy (CMA-ES) (Auger and Hansen, 2005), Opposition-based Tug of War Optimization (OTWO)(Nguyen et al., 2020), Teaching-Learning-Based Optimization (TLBO) (Rao et al., 2011), Life Choice-Based Optimizer (LCBO) (Khatri et al., 2020), Global Mine Blast Algorithm (GMBA) (Yadav et al., 2020), Equilibrium Optimizer (EO) (Faramarzi et al., 2020), Moth-Flame Optimization algorithm (MFO)(Mirjalili, 2015b), Snake Optimizer (SO) (Hashim and Hussien, 2022), Grey Wolf Optimizer (GWO) (Mirjalili et al., 2014), Whale Optimization Algorithm (WOA) (Mirjalili and Lewis, 2016), Giant Trevally Optimizer (GTO) (Sadeeq and Abdulazeez, 2022), Dragonfly Algorithm (DA) (Mirjalili, 2016a), Salp Swarm

Algorithm (SSA) (Mirjalili et al., 2017), Enhanced Colliding Bodies Optimization (ECBO) (Kaveh and Ghazaan, 2014), Multi-Verse Optimizer (MVO) (Mirjalili et al., 2016), Water Evaporation Optimization (WEO) (Kaveh and Bakhshpoori, 2016), Monkey King Evolution (MKE) (Meng and Pan, 2016), and Levy flight Jaya algorithm (LJA) (Iacca et al., 2021). Among these, CMA-ES is an evolutionary algorithm. TLBO, LCBO, and LJA belong to human society-inspired heuristics. OTWO, GMBA, EO, ECBO, MVO, and WEO are inspired by physical, chemical, or mathematical principles. Finally, MFO, SO, GWO, WOA, GTO, DA, SSA, and MKE are categorized as swarm intelligence-based heuristics.

All algorithms are implemented in Java and executed using NetBeans IDE 20 on a system equipped with an Intel(R) Core(TM) i7-1065G7 CPU and 16 GB of RAM. The CMA-ES implementation was obtained from the respective authors' websites, while the remaining algorithms were coded by the author. Readers can also find corresponding MATLAB or Python implementations in open-source platforms like MEALPY (Van Thieu and Mirjalili, 2023), ensuring the reproducibility of results presented in this paper. Parameter settings for all algorithms are detailed in Table 2. Each algorithm uses a population size of 50 and a fixed maximum of 2000 iterations, although CMA-ES often terminates earlier due to its internal convergence mechanism. Other algorithm-specific parameters were set based on the recommendations in their original papers. Each algorithm was executed 30 times, and the mean, standard deviation, and best solution from these independent runs were recorded.

Table 2 The 20 algorithms to be compared and their parameters.

| Algorithm | Parameters |
|---|---|
| CMA-ES | $N$=50, $\mu = 30$, $maxNumIter$=2000, others are set with default values given in paper (Beyer and Sendhoff, 2008) |
| OTWO | $N$=50, $maxNumIter$=2000, static coefficient of friction=1.0, kinematic coefficient of friction=1.2, $\alpha = 0.95$, $\beta = 0.02$ |
| TLBO | $N$=50, $maxNumIter$=2000 |
| LCBO | $N$=50, $maxNumIter$=2000, other parameters are set as the authors suggested in paper (Rao et al., 2011) |
| GMBA | $N$=50, $maxNumIter$=2000, the number of explorations equals 1000, the reduction constant is 120 |
| EO | $N$=50, $maxNumIter$=2000, the size of memory is 4 |
| MFO | $N$=50, $maxNumIter$=2000, the parameter of logarithmic spiral is 1.0 |
| SO | $N$=50, $maxNumIter$=2000. |
| VLA | $N$=50, $maxNumIter$=2000, $CR$=0.25, $n_R$=100. |
| AVLA | $N$=50, $maxNumIter$=2000, $H$=50, $n_R$=100. |
| GWO | $N$=50, $maxNumIter$=2000 |
| WOA | $N$=50, $maxNumIter$=2000, parameter of Logarithmic spiral is 1.1 |
| GTO | $N$=50, $maxNumIter$=2000, other parameters are set as the authors suggested in paper (Sadeeq and Abdulazeez, 2022) |
| DA | $N$=50, $maxNumIter$=2000, Initial Inertia Weight =0.9, Final Inertia Weight=0.4, Initial Enemy Distraction Weight=0.1, final Enemy Distraction Weight=0.0 |
| SSA | $N$=50, $maxNumIter$=2000, other parameters are set as the authors suggested in paper (Mirjalili et al., 2017) |
| ECBO | $N$=50, $maxNumIter$=2000, size of memory=5, mutation probability=0.1 |
| MVO | $N$=50, $maxNumIter$=2000, minimal wormhole existence probability (WEP)=0.2, maximal WEP=1.0, the exploitation accuracy p in travelling distance rate (TDR)=6 |
| WEO | $N$=50, $maxNumIter$=2000, other parameters are set as the authors suggested in paper (Kaveh and Bakhshpoori, |
| | 2016) |
| MKE | *N*=50, *maxNumIter*=2000, FC=0.7 |
| LJA | *N*=50, *maxNumIter*=2000, beta of Levy flight=1.5 |

## 5.2. Test on 29 test benchmarks

The 29 well-known optimization problems considered in this study are grouped into four categories, detailed in Appendix A. The first category comprises seven unimodal problems. The second category consists of six multimodal problems. The third category includes ten fixed-dimension multimodal problems, and the final category contains six composite benchmark functions.

In the subsequent tables, abbreviations are used for clarity: "Ave" denotes the average value (mean), "Std" represents the standard deviation, and "Best" indicates the best objective value obtained from 30 independent runs of each algorithm. "FM" stands for Friedman mean, while "FMR" represents the Friedman mean rank. The values in the "FM of Aves" row are computed as follows: first, the means obtained by the 20 algorithms for each optimization problem are ranked, with the smallest mean ranked first and the largest ranked last. For each algorithm, the ranks across all problems are summed, and the total is divided by the number of ranks. These divided values constitute the "FM of Aves." The ranks in "FMR of Aves" are derived by comparing the relative "FM of Aves" values across the 20 algorithms. Similarly, "FM of Bests" and "FMR of Bests" are calculated based on the best objective values achieved.

The results for the first category of seven unimodal problems are presented in Table 3. According to these data, LCBO, GTO, and SO rank first, second, and third, respectively, in the FMR of Aves. AVLA achieves sixth place in the Friedman mean rank of average objective function values. However, when considering the best solutions obtained across 30 runs, AVLA secures the top position in the FMR of Bests. A higher rank in the FMR of Aves for these unimodal problems indicates a stronger capacity for exploitation. However, excessive exploitation tendencies may result in weaker performance on other problem types. This observation is supported by the subsequent experiments.

The results for the six multimodal problems in the second category are presented in Table 4. Here, AVLA achieves first place in both Friedman mean ranks (Aves and Bests). VLA ranks second in the FMR of Aves and third in the FMR of Bests. LCBO retains its strong performance, securing second place in the FMR of Bests, consistent with its ranking in the first category. WEO emerges in this category, achieving third place in the FMR of Aves.

Results for the ten fixed-dimension multimodal problems in the third category are shown in Table 5. AVLA maintains its outstanding performance, achieving first place in both the FMR of Aves and FMR of Bests. WEO and VLA rank second and third, respectively, in the FMR of Aves. Notably, seven algorithms share the top position in the FMR of Bests, highlighting the competitive nature of performance in this category.

The results for the six composite benchmark problems in the final category are presented in Table 6. WEO, AVLA, and GMBA rank first, second, and third, respectively, in the FMR of Aves, with VLA following closely in fourth place. Interestingly, VLA outperforms AVLA in some cases, such as for problems F24 and F26, suggesting that parameter adaptation does not universally enhance performance. In terms of the FMR of Bests, AVLA, ECBO, and LCBO take the top three positions, while VLA secures fourth place.

To evaluate the overall performance of all 20 algorithms across the 29 problems, we calculated the Friedman mean ranks of Aves and Bests and summarized the results in Table 7. AVLA demonstrates superior performance, ranking first in both the FMR of Aves and the FMR of Bests. VLA also performs well, achieving second place in the FMR of Aves and third place in the FMR of Bests. LCBO ranks second in the FMR of Bests, while WEO secures third place in the FMR of Aves.

In summary, AVLA consistently outperforms the other 19 algorithms based on both Friedman mean ranks, indicating its strong overall effectiveness. VLA also proves to be a competitive algorithm, occasionally providing better solutions than AVLA. Notably, AVLA achieves a well-balanced trade-off between exploration and exploitation, contributing to its robust performance.

Table 3 Optimization results and comparison for unimodal benchmark functions.

| Fun. | Index | CMA-ES | OTWO | TLBO | LCBO | GMBA | EO | MFO | SO | VLA | AVLA |
|---|---|---|---|---|---|---|---|---|---|---|---|
| F1 | Ave | 1.13E04 | 1.89E-01 | 0 | 0 | 1.29E-03 | 1.65E-278 | 2.04E-69 | 0 | 4.4E-33 | 5.71E-84 |
|  | Std | 9E03 | 1.52E-01 | 0 | 0 | 5.92E-03 | 0 | 1.12-68 | 0 | 8.68E-32 | 3.12E-83 |
|  | Best | 1.07E-17 | 1.28E-02 | 0 | 0 | 6.36E-13 | 4.95E-289 | 2.27E-109 | 0 | 6.02E-36 | 0 |
| F2 | Ave | 2.23E01 | 9.38E-02 | 3.66E-290 | 0 | 1.6E-01 | 8.41E-153 | 3.33E-01 | 1.04E-315 | 2.75E-26 | 1.55E-28 |
|  | Std | 1.38E01 | 3.51E-02 | 0E00 | 0 | 5.56E-01 | 2.1E-152 | 1.83E00 | 0E00 | 4.13E-26 | 1.57E-28 |
|  | Best | 7.36E-10 | 3.33E-02 | 1.07E-303 | 0 | 2.93E-07 | 9.86E-156 | 4.42E-54 | 4.9E-324 | 1.39E-27 | 0 |
| F3 | Ave | 4.71E05 | 4.58E00 | 0 | 0 | 1.77E-03 | 2.39E-278 | 3.33E02 | 0 | 2.17E-30 | 0 |
|  | Std | 3.31E05 | 2.79E00 | 0 | 0 | 4.98E-03 | 0E00 | 1.83E03 | 0 | 5.49E-30 | 0 |
|  | Best | 1.85E-17 | 1.36E00 | 0 | 0 | 8.25E-09 | 1E-284 | 3.7E-109 | 0 | 1.3E-33 | 0 |
| F4 | Ave | 1E02 | 2.03E-01 | 1.53E-241 | 0 | 9.37E-04 | 2.61E-98 | 1.86E-09 | 4.13E-309 | 1.1E-01 | 9.16E-22 |
|  | Std | 0E00 | 8.49E-02 | 0E00 | 0 | 2.28E-03 | 1.41E-97 | 4.27E-09 | 0E00 | 9.46E-02 | 3.82E-21 |
|  | Best | 1E02 | 6.63E-02 | 1.51E-254 | 0 | 6.33E-07 | 7.34E-106 | 1.69E-13 | 2.36E-314 | 8.25E-03 | 0 |
| F5 | Ave | 8.97E07 | 9.05E00 | 6.69E00 | 3.9E00 | 3.28E01 | 4.97E00 | 9.15E03 | 8.87E00 | 1.81E00 | 8.65E-29 |
|  | Std | 7.16E07 | 3.22E00 | 1.02E00 | 3.44E-01 | 6.91E01 | 8.71E-01 | 2.74E04 | 9.04E-02 | 1.82E00 | 3.58E-28 |
|  | Best | 5.5E-19 | 2.67E00 | 3.3E00 | 3.12E00 | 2.52E-01 | 3.7E00 | 5.08E-02 | 8.64E00 | 2.23E-03 | 0 |
| F6 | Ave | 4.33E03 | 6.67E-02 | 0 | 0 | 1.2E00 | 0 | 3E-01 | 0 | 0 | 0 |
|  | Std | 5.68E03 | 3.65E-01 | 0 | 0 | 3.96E00 | 0 | 7.94E-01 | 0 | 0 | 0 |
|  | Best | 0 | 0 | 0 | 0 | 0 | 0 | 0E00 | 0 | 0 | 0 |
| F7 | Ave | 1.01E01 | 4.3E-05 | 1.84E-04 | 3.29E-05 | 1.14E-02 | 5.97E-05 | 4.49E-03 | 2.05E-05 | 3.45E-05 | 3.46E-04 |
|  | Std | 1.04E01 | 3.44E-05 | 1.38E-04 | 1.99E-05 | 3.03E-02 | 4.68E-05 | 4.1E-03 | 1.72E-05 | 2.35E-05 | 1.53E-04 |
|  | Best | 7.7E-03 | 5.52E-06 | 2.69E-05 | 2.7E-06 | 2.73E-04 | 4.73E-06 | 5.29E-04 | 4.15E-07 | 2.48E-06 | 9.72E-05 |
| FM of Aves |  | 0.095238 | 0.065306 | 0.02381 | 0.015986 | 0.070748 | 0.026871 | 0.065986 | 0.020408 | 0.039796 | 0.028912 |
| FMR of Aves |  | 20 | 13 | 4 | 1 | 17 | 5 | 14 | 3 | 8 | 6 |
| FM of Bests |  | 0.060544 | 0.069388 | 0.030612 | 0.02449 | 0.059864 | 0.036054 | 0.044898 | 0.031293 | 0.042177 | 0.020408 |
| FMR of Bests |  | 15 | 18 | 4 | 2 | 14 | 6 | 9 | 5 | 7 | 1 |

Continue the above Table 3:

| Fun. | Index | GWO | WOA | GTO | DA | SSA | ECBO | MVO | WEO | MKE | LJA |
|---|---|---|---|---|---|---|---|---|---|---|---|
| F1 | Ave | 1.62E-143 | 6.3E-98 | 0 | 2.09E03 | 1.47E-09 | 5.82E-20 | 6.02E-04 | 1.88E-18 | 1.01E-11 | 7E-08 |
|  | Std | 3.89E-143 | 2.21E-97 | 0 | 3.28E03 | 9.77E-10 | 1.13E-19 | 2.71E-04 | 2.82E-18 | 3.42E-11 | 5.22E-08 |
|  | Best | 1.88E-148 | 7.08E-108 | 0 | 1.01E01 | 3.87E-10 | 1.4E-26 | 2.35E-04 | 9.48E-20 | 2.24E-14 | 8.28E-09 |
| F2 | Ave | 8.21E-77 | 2.9E-57 | 0 | 2.2E01 | 7.3E-02 | 5.3E-12 | 6.42E-03 | 2.89E-07 | 2.67E00 | 6E-05 |
|  | Std | 1.45E-76 | 9.78E-57 | 0 | 5.27E00 | 1.24E-01 | 5.02E-12 | 1.27E-03 | 2.53E-07 | 5.21E00 | 4.25E-05 |
|  | Best | 4.06E-79 | 6.2E-61 | 0 | 9.63E00 | 2.32E-04 | 1.2E-13 | 3.4E-03 | 5.11E-08 | 1.31E-07 | 1.93E-05 |
| F3 | Ave | 6.35E-141 | 2.39E-97 | 0 | 6.46E04 | 9.05E01 | 4.36E-20 | 5.87E-01 | 4.44E-17 | 1.37E04 | 1.19E-06 |
|  | Std | 2.5E-14. | 7.21E-97 | 0 | 7.5E04 | 1.14E02 | 1.36E-19 | 7.35E-01 | 4.64E-17 | 3.08E04 | 2.13E-06 |
|  | Best | 9.76E-148 | 1.86E-105 | 0 | 1.85E03 | 5.26E-02 | 5.45E-28 | 3.07E-02 | 1.73E-18 | 1.42E-12 | 5.94E-08 |
| F4 | Ave | 1.13E-48 | 2.1E-30 | 0 | 2.93E01 | 4.84E-02 | 6.45E-11 | 1.55E-02 | 6.21E-04 | 2.22E-05 | 2.5E-02 |
|  | Std | 3.74E-48 | 7.02E-30 | 0 | 1.43E01 | 9.37-02 | 2.46E-10 | 6.53E-03 | 2.44E-04 | 3.82E-05 | 1.05E-02 |
|  | Best | 3.59E-52 | 2.7E-34 | 0 | 7.38E00 | 2.14E-05 | 1.15E-013 | 7.29E-03 | 1.55E-04 | 4.41E-07 | 9.77E-03 |
| F5 | Ave | 7.11E00 | 5.19E00 | 7.85E00 | 1.02E06 | 7.27E01 | 3.62E00 | 7.89E01 | 1.62E00 | 1.23E04 | 4.13E01 |
|  | Std | 2.96E-01 | 7.06E-01 | 1.37E-01 | 2.78E06 | 1.89E02 | 3.47E00 | 1.46E02 | 1.04E00 | 3.1E04 | 4.8E01 |
|  | Best | 6.24E00 | 3.86E00 | 7.63E00 | 2.54E02 | 2.43E-03 | 1.7E-02 | 3.53E-01 | 2.98E-01 | 1.32E-02 | 4.74E00 |
| F6 | Ave | 0 | 0 | 0 | 2.1E03 | 1.8E00 | 0 | 4E-01 | 0 | 3.8E00 | 0 |
|  | Std | 0 | 0 | 0 | 3.48E03 | 1.32E00 | 0 | 6.21E-01 | 0 | 6.69E00 | 0 |
|  | Best | 0 | 0 | 0 | 1.4E01 | 0E00 | 0 | 0E00 | 0 | 0E00 | 0 |
| F7 | Ave | 1.52E-02 | 2.24E-02 | 1.03E-05 | 6.25E-01 | 1.28E-03 | 4.83E-04 | 5.86E-04 | 3.92E-03 | 6.54E-03 | 7.95E-03 |
|  | Std | 1.2E-02 | 2.14E-02 | 1.06E-05 | 6.79E-01 | 8.82E-04 | 3.08E-04 | 3.45E-04 | 1.69E-03 | 6.89E-03 | 4.23E-03 |
|  | Best | 1.45E-04 | 1.4E-03 | 3.88E-07 | 2.67E-02 | 1.87E-04 | 9.59E-05 | 1.23E-04 | 1.46E-03 | 6.07E-04 | 1.94E-03 |
| FM of Aves |  | 0.039116 | 0.041156 | 0.018027 | 0.090476 | 0.070068 | 0.040476 | 0.066667 | 0.045238 | 0.07483 | 0.060884 |
| FMR of Aves |  | 7 | 10 | 2 | 19 | 16 | 9 | 15 | 11 | 18 | 12 |
| FM of Bests |  | 0.044898 | 0.05102 | 0.026531 | 0.094558 | 0.061224 | 0.044898 | 0.066667 | 0.059184 | 0.056463 | 0.07483 |
| FMR of Bests |  | 9 | 11 | 3 | 20 | 16 | 9 | 17 | 13 | 12 | 19 |

Table 4 Optimization results and comparison for multimodal benchmark functions.

| Fun. | Index | CMA-ES | OTWO | TLBO | LCBO | GMBA | EO | MFO | SO | VLA | AVLA |
|---|---|---|---|---|---|---|---|---|---|---|---|
| F8 | Ave | -2.05E03 | -3.23E03 | -3.24E03 | -3.59E03 | -2.08E03 | -3.25E03 | -3.04E03 | -2.55E03 | -4.06E03 | -4.19E03 |
|  | Std | 1.04E02 | 4.55E02 | 2.35E02 | 2.66E02 | 3.84E02 | 2.43E02 | 3.28E02 | 2.02E02 | 6.29E01 | 2.78E-12 |
|  | Best | -2.49E03 | -4.19E03 | -3.83E03 | -3.97E03 | -3.48E03 | -3.71E03 | -3.6E03 | -2.99E03 | -4.19E03 | -4.19E03 |
| F9 | Ave | 6.65E01 | 1.42E-01 | 1.04E00 | 5.97E-01 | 2.04E01 | 1.67E-01 | 1.7E01 | 0 | 1.24E-02 | 0 |
|  | Std | 3.55E01 | 3.65E-01 | 1.82E00 | 2.15E00 | 1.02E01 | 9.15E-01 | 1.18E01 | 0 | 6.68E-02 | 0 |
|  | Best | 1.39E01 | 8.36E-03 | 0 | 0 | 9.95E00 | 0 | 3.98E00 | 0 | 0 | 0 |
| F10 | Ave | 2E01 | 3.73E-01 | 3.88E-15 | 4.44E-16 | 1.46E01 | 3.29E-15 | 1.25E-14 | 4.44E-16 | 2.12E-14 | 4.44E-16 |

|  | Std | 5.01E-03 | 2.69E-01 | 6.49E-16 | 0 | 8.98E00 | 1.45E-15 | 7.95E-15 | 0E00 | 2.32E-14 | 0E00 |
|  | Best | 1.99E01 | 6.67E-02 | 4.44E-16 | 4.44E-16 | 3.98E-07 | 4.44E-16 | 4E-15 | 4.44E-16 | 4E-15 | 4.44E-16 |
| F11 | Ave | 6.38E01 | 2.07E-01 | 1.33E-02 | 3.28E-03 | 2.47E-01 | 1.65E-02 | 1.91E-01 | 0 | 1.32E-03 | 1.58E-52 |
|  | Std | 5.6E)1 | 1.86E-01 | 2.33E-02 | 1.18E-02 | 2.37E-01 | 1.64E-02 | 1.33E-01 | 0 | 2.78E-03 | 8.15E-52 |
|  | Best | 1.09E-19 | 1.25E-02 | 0 | 0 | 1.3E-286 | 6.4E-02 | 3.2E-02 | 0 | 2.15E-08 | 0 |
| F12 | Ave | 2.29E08 | 5.78E-01 | 1.33E-04 | 3.11E-02 | 9.9E-07 | 5E-20 | 1.55E-01 | 4.51E-01 | 5.63E-31 | 4.71E-32 |
|  | Std | 2.73E08 | 9.46E-01 | 5.76E-04 | 1.25E-01 | 3.08E-06 | 1.15E-19 | 3.54E-01 | 1.51E-01 | 2E-30 | 1.67E-47 |
|  | Best | 4.18E-20 | 3.52E-08 | 8.97E-14 | 4.71E-32 | 1.51E-14 | 3.91E-23 | 4.71E-32 | 1.91E-01 | 4.71E-32 | 4.71E-32 |
| F13 | Ave | 3.11E08 | 1.05E-01 | 4.86E-02 | 1.52E-01 | 3.9E-03 | 8.08E-03 | 8.37E-02 | 9.17E-01 | 2.76E-24 | 1.35E-32 |
|  | Std | 3.13E08 | 1.26E-01 | 3.02E-02 | 1.1E-01 | 2.07E-02 | 3.08E-02 | 4.14E-01 | 1.48E-01 | 1.51E-23 | 5.57E-48 |
|  | Best | 2.81E-19 | 1.91E-07 | 3.12E-12 | 1.35E-32 | 4.08E-14 | 2.14E-21 | 1.56E-32 | 4.84E-01 | 2.19E-32 | 1.35E-32 |
| FM of Aves |  | 0.093651 | 0.057937 | 0.042857 | 0.04127 | 0.062698 | 0.034127 | 0.056349 | 0.045635 | 0.020635 | 0.009524 |
| FMR of Aves |  | 20 | 14 | 8 | 7 | 17 | 4 | 13 | 10 | 2 | 1 |
| FM of Bests |  | 0.062698 | 0.058333 | 0.03373 | 0.017857 | 0.063889 | 0.030556 | 0.044048 | 0.054365 | 0.02619 | 0.015079 |
| FMR of Bests |  | 15 | 14 | 6 | 2 | 16 | 4 | 7 | 12 | 3 | 1 |

Continue the above Table 4:

| Fun. | Index | GWO | WOA | GTO | DA | SSA | ECBO | MVO | WEO | MKE | LJA |
|---|---|---|---|---|---|---|---|---|---|---|---|
| F8 | Ave | -8.6E02 | -2.97E03 | -2.72E03 | -2.26E03 | -2.87E03 | -3.16E03 | -3.02E03 | -4.17E03 | -3.27E03 | -2.86E03 |
|  | Std | 1.29E02 | 2.97E02 | 4.72E02 | 3.4E02 | 3.66E02 | 5.44E02 | 3.16E02 | 3.21E01 | 3.17E02 | 1.97E02 |
|  | Best | -1.24E03 | -3.74E03 | -3.7E03 | -3.26E03 | -3.48E03 | -4.07E03 | -3.26E03 | -4.19E03 | -3.83E03 | -3.24E03 |
| F9 | Ave | 0 | 0 | 0 | 7.68E01 | 9.32E00 | 1.88E01 | 1.71E01 | 1.62E-01 | 3.49E01 | 4.12E01 |
|  | Std | 0 | 0 | 0 | 9.56E00 | 4.23E00 | 7.96E00 | 1.23E01 | 1.91E-01 | 1.32E01 | 7.81E00 |
|  | Best | 0 | 0 | 0 | 5.23E01 | 3.98E00 | 0E00 | 2.99E00 | 8.38E-03 | 1.79E01 | 2.2E01 |
| F10 | Ave | 4.44E-16 | 1.93E01 | 4.44E-16 | 1.84E01 | 1.29E00 | 2.65E-12 | 1.19E-01 | 1.16E-05 | 9.01E-01 | 1.68E01 |
|  | Std | 0E00 | 3.64E00 | 0E00 | 7.78E-01 | 9.52E-01 | 3.3E-12 | 4.15E-01 | 1.38E-05 | 1.23E00 | 2.34E00 |
|  | Best | 4.44E-16 | 3.92E-12 | 4.44E-16 | 1.64E01 | 8.62E-06 | 1.53E-13 | 4.26E-03 | 5.48E-07 | 3.41E-08 | 8.81E00 |
| F11 | Ave | 1.64E-03 | 4.84E-03 | 0 | 1.07E01 | 2.96E-01 | 3.77E-01 | 3.18E-01 | 3.43E-02 | 2.03E-01 | 5.67E-01 |
|  | Std | 7.02E-03 | 9.31E-03 | 0 | 1.51E01 | 1.6E-01 | 1.47E-01 | 1.17E-01 | 1.42E-02 | 9.58E-02 | 1.15E-01 |
|  | Best | 6.89E-139 | 5.56E-99 | 0 | 7.81E-01 | 4.43E-02 | 9.87E-05 | 7.47E-02 | 1.22E-02 | 4.68E-02 | 2.85E-01 |
| F12 | Ave | 1.11E-02 | 3.27E-03 | 8.74E-02 | 3.16E06 | 1.39E00 | 2.78E-17 | 1.04E-02 | 1.98E-13 | 8.29E-02 | 8.19E-07 |
|  | Std | 1.43E-02 | 9.04E-03 | 3.47E-02 | 7.34E06 | 1.56E00 | 1.45E-16 | 5.68E-02 | 3.69E-13 | 2.57E-01 | 2.51E-06 |
|  | Best | 1.61E-08 | 3.49E-07 | 3.28E-02 | 1.88E00 | 6.83E-09 | 1.15E-23 | 5.75E-06 | 1.39E-15 | 4.3E-16 | 1.75E-09 |
| F13 | Ave | 4.67E-03 | 7.5E-06 | 2.65E-02 | 1.39E06 | 8.83E-03 | 6.58E-19 | 1.29E-04 | 1.14E-14 | 1.5E-01 | 5.88E-07 |
|  | Std | 2.23E-02 | 4.1E-09 | 7.22E-03 | 7.31E06 | 1.66E-02 | 2.9E-18 | 8.62E-05 | 1.74E-14 | 4.65E-01 | 9.65E-07 |
|  | Best | 1.32E-07 | 2.37E-06 | 5.43E-03 | 3.35E00 | 2.26E-09 | 2.64E-24 | 2.31E-05 | 2.13E-17 | 8.49E-15 | 4.38E-08 |
| FM of Aves |  | 0.040476 | 0.044444 | 0.038492 | 0.088889 | 0.066667 | 0.046032 | 0.055556 | 0.031746 | 0.061111 | 0.061905 |
| FMR of Aves |  | 6 | 9 | 5 | 19 | 18 | 11 | 12 | 3 | 15 | 16 |
| FM of Bests |  | 0.05119 | 0.052778 | 0.048016 | 0.09127 | 0.066667 | 0.032937 | 0.070635 | 0.045635 | 0.056349 | 0.077778 |
| FMR of Bests |  | 10 | 11 | 9 | 20 | 17 | 5 | 18 | 8 | 13 | 19 |

**Table 5** Optimization results for fixed dimension multimodal benchmark functions.

| Fun. | Index | CMA-ES | OTWO | TLBO | LCBO | GMBA | EO | MFO | SO | VLA | AVLA |
|---|---|---|---|---|---|---|---|---|---|---|---|
| F14 | Ave | 2.58E01 | 3.33E00 | 1.03E00 | 2.37E00 | 9.98E-01 | 3.35E00 | 1.53E00 | 1.61E00 | 9.98E-01 | 9.98E-01 |
|  | Std | 8.97E01 | 2.21E00 | 1.81E-01 | 3.37E00 | 1.81E-08 | 3.67E00 | 1.24E00 | 7.83E-01 | 1.13E-16 | 1.13E-16 |
|  | Best | 9.98E-01 | 9.98E-01 | 9.98E-01 | 9.98E-01 | 9.98E-01 | 9.98E-01 | 9.98E-01 | 9.98E-01 | 9.98E-01 | 9.98E-01 |
| F15 | Ave | 1.44E-02 | 2.52E-03 | 4.42E-04 | 2.34E-03 | 7.06E-04 | 1.68E-03 | 1.24E-03 | 6.72E-04 | 3.51E-04 | 3.09E-04 |
|  | Std | 8.18E-03 | 6.06E-03 | 3E-04 | 6.11E-03 | 4.13E-04 | 5.08E-03 | 1.43E-03 | 2.36E-04 | 4.61E-05 | 4.07E-06 |
|  | Best | 2.25E-03 | 3.07E-04 | 3.07E-04 | 3.07E-04 | 3.07E-04 | 3.07E-04 | 3.08E-04 | 3.37E-04 | 3.08E-04 | 3.07E-04 |
| F16 | Ave | 1.39E01 | -1.03E00 | -1.03E00 | -1.03E00 | -1.03E00 | -1.03E00 | -1.03E00 | -1.03E00 | -1.03E00 | -1.03E00 |
|  | Std | 8.21E01 | 6.9E-09 | 5.86E-16 | 6.71E-16 | 2.95E-11 | 6.78E-16 | 6.25E-16 | 2.31E-04 | 6.78E-16 | 6.78E-16 |
|  | Best | -1.03E00 | -1.03E00 | -1.03E00 | -1.03E00 | -1.03E00 | -1.03E00 | -1.03E00 | -1.03E00 | -1.03E00 | -1.03E00 |
| F17 | Ave | 1.44E00 | 3.98E-01 | 3.98E-01 | 3.98E-01 | 3.98E-01 | 3.98E-01 | 3.98E-01 | 4.04E-01 | 3.98E-01 | 3.98E-01 |
|  | Std | 2.12E00 | 1.15E-09 | 0E00 | 0E00 | 1.66E-08 | 0E00 | 0E00 | 7.8E-03 | 0E00 | 0E00 |
|  | Best | 3.98E-01 | 3.98E-01 | 3.98E-01 | 3.98E-01 | 3.98E-01 | 3.98E-01 | 3.98E-01 | 3.98E-01 | 3.98E-01 | 3.98E-01 |
| F18 | Ave | 4.37E01 | 3E00 | 3E00 | 3E00 | 1.11E01 | 3E00 | 3E00 | 3.03E00 | 3E00 | 3E00 |
|  | Std | 7.35E01 | 5.58E-06 | 2.26E-15 | 2.44E-15 | 1.26E01 | 3.38E-09 | 3.02E-15 | 3.34E-02 | 2.26E-15 | 2.26E-15 |
|  | Best | 3E00 | 3E00 | 3E00 | 3E00 | 3E00 | 3E00 | 3E00 | 3E00 | 3E00 | 3E00 |
| F19 | Ave | -3.19E00 | -3.78E00 | -3.78E00 | -3.78E00 | -3.54E00 | -3.78E00 | -3.78E00 | -3.76E00 | -3.78E00 | -3.78E00 |
|  | Std | 8.24E-01 | 2.53E-15 | 2.26E-15 | 2.26E-15 | 6.17E-01 | 2.26E-15 | 2.26E-15 | 9.2E-03 | 2.26E-15 | 2.26E-15 |
|  | Best | -3.78E00 | -3.78E00 | -3.78E00 | -3.78E00 | -3.78E00 | -3.78E00 | -3.78E00 | -3.78E00 | -3.78E00 | -3.78E00 |
| F20 | Ave | -5.14E-01 | -5.72E-01 | -5.73E-01 | -5.73E-01 | -5.72E-01 | -5.73E-01 | -5.73E-01 | -5.11E-01 | -5.73E-01 | -5.73E-01 |
|  | Std | 1.11E-01 | 1.52E-04 | 3.39E-16 | 3.57E-16 | 2.63E-04 | 3.22E-12 | 2.14E-16 | 1.45E-02 | 3.43E-16 | 2.4E-16 |
|  | Best | -5.73E-01 | -5.73E-01 | -5.73E-01 | -5.73E-01 | -5.73E-01 | -5.73E-01 | -5.73E-01 | -5.59E-01 | -5.73E-01 | -5.73E-01 |
| F21 | Ave | -4.97E00 | -8.48E00 | -8.86E00 | -6.58E00 | -7.05E00 | -1.02E01 | -6.65E00 | -3.47E00 | -1.02E01 | -1.02E01 |
|  | Std | 4.17E00 | 3.11E00 | 2.45E00 | 3.68E00 | 3.27E00 | 5.42E-15 | 3.45E00 | 1.76E00 | 8.42E-15 | 8.76E-15 |
|  | Best | -1.02E01 | -1.02E01 | -1.02E01 | -1.02E01 | -1.02E01 | -1.02E01 | -1.02E01 | -8.12E00 | -1.02E01 | -1.02E01 |
| F22 | Ave | -4.49E00 | -8.46E00 | -9.6E00 | -7.15E00 | -6.6E00 | -3.7E00 | -8.38E00 | -3.83E00 | -1.04E01 | -1.04E01 |
|  | Std | 4.15E00 | 3.3E00 | 1.8E00 | 3.79E00 | 3.5E00 | 3.08E-12 | 3.45E00 | 1.5E00 | 0E00 | 0E00 |
|  | Best | -1.04E01 | -1.04E01 | -1.04E01 | -1.04E01 | -1.04E01 | -3.7E00 | -1.04E01 | -7.61E00 | -1.04E01 | -1.04E01 |
| F23 | Ave | -3.85E00 | -8.42E00 | -9.79E00 | -7E00 | -4.91E00 | -1.05E01 | -7.61E00 | -3.65E00 | -1.05E01 | -1.05E01 |
|  | Std | 3.51E00 | 3.06E00 | 1.98E00 | 3.42E00 | 2.99E00 | 5.42E-15 | 3.45E00 | 1.52E00 | 3.57E-15 | 3.28E-15 |
|  | Best | -1.05E01 | -1.05E01 | -1.05E01 | -1.05E01 | -1.05E01 | -1.05E01 | -1.05E01 | -8.89E00 | -1.05E01 | -1.05E01 |
| FM of Aves |  | 0.091429 | 0.051905 | 0.035476 | 0.055238 | 0.062143 | 0.045238 | 0.049524 | 0.075238 | 0.028095 | 0.026667 |
| FMR of Aves |  | 20 | 13.5 | 4 | 15 | 16 | 9 | 10.5 | 18 | 3 | 1 |
| FM of Bests |  | 0.0519 | 0.0448 | 0.0448 | 0.0448 | 0.0448 | 0.0498 | 0.0481 | 0.0686 | 0.0481 | 0.0448 |
| FMR of Bests |  | 16 | 4 | 4 | 4 | 4 | 13 | 10 | 20 | 10 | 4 |

Continue the above Table 5:

| Fun. | Index | GWO | WOA | GTO | DA | SSA | ECBO | MVO | WEO | MKE | LJA |
|---|---|---|---|---|---|---|---|---|---|---|---|
| F14 | Ave | 8.22E00 | 3.36E00 | 9.98E-01 | 6.14E00 | 1.03E00 | 9.98E-01 | 9.98E-01 | 9.98E-01 | 1.06E00 | 9.98E-01 |
| | Std | 4.45E00 | 2.63E00 | 3E-05 | 4.44E00 | 1.81E-01 | 2.48E-10 | 3.3E-10 | 1.01E-16 | 3.62E-01 | 4.12E-17 |
| | Best | 2.98E00 | 9.98E-01 | 9.98E-01 | 9.98E-01 | 9.98E-01 | 9.98E-01 | 9.98E-01 | 9.98E-01 | 9.98E-01 | 9.98E-01 |
| F15 | Ave | 1.25E-02 | 1.04E-03 | 3.17E-04 | 2.4E-03 | 1.41E-03 | 1.45E-03 | 4.67E-03 | 3.43E-04 | 4.71E-03 | 2.89E-03 |
| | Std | 1.01E-02 | 3.76E-04 | 7.89E-06 | 1.3E-03 | 3.68E-03 | 3.58E-03 | 7.98E-03 | 4.1E-05 | 1.15E-02 | 4.33E-03 |
| | Best | 3.07E-04 | 3.07E-04 | 3.1E-04 | 9.14E-04 | 3.08E-04 | 4.79E-04 | 3.08E-04 | 3.08E-04 | 3.07E-04 | 3.43E-04 |
| F16 | Ave | -1.03E00 | -1.03E00 | -1.03E00 | -1.03E00 | -1.03E00 | -1.03E00 | -1.03E00 | -1.03E00 | -1.03E00 | -1.03E00 |
| | Std | 4.51E-09 | 5.08E-09 | 1.09E-06 | 6.63E-03 | 2.17E-15 | 1.33E-04 | 2.05E-05 | 6.58E-16 | 6.25E-16 | 8.67E-07 |
| | Best | -1.03E00 | -1.03E00 | -1.03E00 | -1.03E00 | -1.03E00 | -1.03E00 | -1.03E00 | -1.03E00 | -1.03E00 | -1.03E00 |
| F17 | Ave | 3.98E-01 | 3.98E-01 | 3.98E-01 | 3.98E-01 | 3.98E-01 | 3.98E-01 | .98E-01 | 3.98E-01 | 3.98E-01 | 3.98E-01 |
| | Std | 7.89E-08 | 2.06E-07 | 4.93E-05 | 1.05E-03 | 8.6E-16 | 0E00 | 2.94E-05 | 0E00 | 0E00 | 0E00 |
| | Best | 3.98E-01 | 3.98E-01 | 3.98E-01 | 3.98E-01 | 3.98E-01 | 3.98E-01 | 3.98E-01 | 3.98E-01 | 3.98E-01 | 3.98E-01 |
| F18 | Ave | 3E00 | 3E00 | 3E00 | 3E00 | 3.01E00 | 3E00 | 3E00 | 3E00 | 3E00 | 3E00 |
| | Std | 7.94E-07 | 1.81E-06 | 3.46E-07 | 3.62E-02 | 3.19E-14 | 1.82E-15 | 9.27E-05 | 2.26E-15 | 2.61E-15 | 1.84E-15 |
| | Best | 3E00 | 3E00 | 3E00 | 3E00 | 3E00 | 3E00 | 3E00 | 3E00 | 3E00 | 3E00 |
| F19 | Ave | -3.62E00 | -3.78E00 | -3.78E00 | -3.71E00 | -3.78E00 | -3.78E00 | -3.78E00 | -3.78E00 | -3.78E00 | -3.77E00 |
| | Std | 3.51E-01 | 1.72E-02 | 1.72E-05 | 6E-02 | 3.46E-13 | 2.48E-15 | 2.66E-05 | 2.26E-15 | 2.26E-15 | 3.26E-02 |
| | Best | -3.69E00 | -3.78E00 | -3.78E00 | -3.77E00 | -3.78E00 | -3.78E00 | -3.78E00 | -3.78E00 | -3.78E00 | -3.78E00 |
| F20 | Ave | -5.54E-01 | -5.73E-01 | -5.73E-01 | -5.52E-01 | -5.72E-01 | -5.73E-01 | -5.73E-01 | -5.73E-01 | -5.73E-01 | -5.73E-01 |
| | Std | 3.19E-02 | 1.51E-11 | 1.51E-08 | 2.61E-02 | 2.19E-04 | 2.14E-16 | 3.46E-08 | 2.95E-16 | 2.8E-16 | 3.39E-16 |
| | Best | -5.73E-01 | -5.73E-01 | -5.73E-01 | -5.73E-01 | -5.73E-01 | -5.73E-01 | -5.73E-01 | -5.73E-01 | -5.73E-01 | -5.73E-01 |
| F21 | Ave | -7.41E00 | -8.16E00 | -8.39E00 | -6.67E-01 | -8.9E00 | -7.97E00 | -7.36E00 | -1.02E01 | -6.96E00 | -6.43E00 |
| | Std | 3.66E00 | 3.36E00 | 1.22E00 | 4.34E-01 | 2.61E00 | 2.8E00 | 2.69E00 | 6.36E-15 | 3.16E00 | 2.34E00 |
| | Best | -1.02E01 | -1.02E01 | -1.01E01 | -2.25E00 | -1.02E01 | -1.02E01 | -1.02E01 | -1.02E01 | -1.02E01 | -1.02E01 |
| F22 | Ave | -9.51E00 | -7.82E00 | -9.37E00 | -1.4E00 | -9.21E00 | -1.04E01 | -8.99E00 | -1.04E01 | -7.58E00 | -8.5E00 |
| | Std | 2.32E00 | 3.47E00 | 1.03E00 | 1.81E00 | 2.75E00 | 2.52E-07 | 2.38E00 | 3.93E-15 | 3.59E00 | 2.7E0 |
| | Best | -1.04E01 | -1.04E01 | -1.04E01 | -1.04E01 | -1.04E01 | -1.04E01 | -1.04E01 | -1.04E01 | -1.04E01 | -1.04E01 |
| F23 | Ave | -6.66E00 | -7.26E00 | -8.33E00 | -1.17E00 | -8.69E00 | -7.92E00 | -9.95E00 | -1.05E01 | -8.63E00 | -9.32E00 |
| | Std | 3.22E00 | 3.34E00 | 2.45E00 | 7.75E-01 | 3.13E00 | 2.85E00 | 1.65E00 | 5.73E-14 | 3.24E00 | 2.28E00 |
| | Best | -1.05E01 | -1.05E01 | -1.05E01 | -4.55E00 | -1.05E01 | -1.05E01 | -1.05E01 | -1.05E01 | -1.05E01 | -1.05E01 |
| FM of Aves | | 0.064048 | 0.049524 | 0.035952 | 0.078095 | 0.043571 | 0.039524 | 0.03881 | 0.027619 | 0.051905 | 0.05 |
| FMR of Aves | | 17 | 10.5 | 5 | 19 | 8 | 7 | 6 | 2 | 13.5 | 12 |
| FM of Bests | | 0.0545 | 0.0448 | 0.0538 | 0.0662 | 0.0481 | 0.0510 | 0.0481 | 0.0481 | 0.0448 | 0.0505 |
| FMR of Bests | | 18 | 4 | 17 | 19 | 10 | 15 | 10 | 10 | 4 | 14 |

**Table 6** Optimization results and comparison for composite benchmark functions.

| Fun. | Index | CMA-ES | OTWO | TLBO | LCBO | GMBA | EO | MFO | SO | VLA | AVLA |
|---|---|---|---|---|---|---|---|---|---|---|---|
| F24 | Ave | 2.66E02 | 8.67E01 | 6.33E01 | 5.85E-03 | 7.72E-05 | 7.34E01 | 3.2E00 | 2.29E02 | 5.34E00 | 6.67E00 |
| | Std | 1.48E02 | 7.76E01 | 6.15E01 | 1.93E-02 | 3.97E-04 | 7.85E01 | 4.39E00 | 5.97E01 | 1.81E01 | 2.54E01 |
| | Best | 3.5E01 | 2.09E-08 | 9.35E-11 | 0 | 3.19E-15 | 3.45E-19 | 0E00 | 1.27E02 | 9.63E-29 | 00E00 |
| F25 | Ave | 2.94E02 | 7.86E01 | 5.78E01 | 4.33E01 | 1.03E01 | 2.8E01 | 4.55E01 | 2.58E02 | 1.91E01 | 0 |
| | Std | 7.75E01 | 9.04E01 | 4.58E01 | 2.27E01 | 8.23E00 | 1.59E01 | 5.12E01 | 5.12E01 | 1.54E01 | 0 |
| | Best | 1.9E02 | 9.9E-03 | 4.34E00 | 5.38E-29 | 2.09E-11 | 8.61E-18 | 5.39E00 | 1.39E02 | 2.14E-07 | 0 |
| F26 | Ave | 1.68E03 | 5.72E02 | 2.27E02 | 1.57E02 | 4.18E02 | 2.38E02 | 5.18E02 | 9E02 | 3.04E00 | 4E01 |
| | Std | 1.55E02 | 6.18E02 | 2.29E02 | 2.24E02 | 2.69E02 | 2.37E02 | 2.05E02 | 0E00 | 1.43E01 | 1.18E02 |
| | Best | 1.23E03 | 6.64E-04 | 1.78E-05 | 0 | 1.07E-12 | 1.93E-18 | 3.33E-30 | 9E02 | 1.7E-20 | 0 |
| F27 | Ave | 6.32E01 | 1.12E01 | 6.31E00 | 6.31E00 | 6.31E00 | 1.5E01 | 6.31E00 | 8.64E01 | 6.31E00 | 6.31E00 |
| | Std | 7.86E01 | 2.38E01 | 3.67E-03 | 1.76E-03 | 1.84E-05 | 3.31E01 | 2.88E-04 | 2.11E01 | 2.46E-04 | 3.64E-04 |
| | Best | 6.31E00 | 6.31E00 | 6.31E00 | 6.31E00 | 6.31E00 | 6.31E00 | 6.31E00 | 5.22E01 | 6.31E00 | 6.31E00 |
| F28 | Ave | 2.69E01 | 3.03E01 | 2.1E01 | 6.39E01 | 5.12E-01 | 6.39E01 | 5.44E00 | 6.76E01 | 1.17E00 | 1.85E-01 |
| | Std | 4.28E01 | 5.12E01 | 4.48E01 | 8.43E-01 | 2.86E-01 | 6.02E-01 | 2.16E01 | 3.47E01 | 3.71E-01 | 1.8E-01 |
| | Best | 4.24E00 | 1.25E00 | 5.66E-01 | 7.06E-01 | 1.76E-01 | 1.46E-01 | 6.28E-01 | 1.49E01 | 5.05E-01 | 0 |
| F29 | Ave | 1.01E03 | 5.27E02 | 5.29E02 | 6.52E02 | 4.54E02 | 4.74E02 | 7.19E02 | 7.77E02 | 2.32E02 | 3.6E02 |
| | Std | 1.07E02 | 1.87E02 | 1.68E02 | 2.91E02 | 3.28E02 | 3.24E02 | 2.2E02 | 5.67E01 | 1.71E02 | 1.92E02 |
| | Best | 8.29E02 | 6.97E01 | 4E02 | 0 | 5.19E-11 | 5.28E00 | 0 | 6.81E02 | 8.82E-18 | 0E00 |
| FM of Aves | | 0.088095 | 0.060317 | 0.04246 | 0.03373 | 0.022619 | 0.05119 | 0.044841 | 0.088095 | 0.023413 | 0.017857 |
| FMR of Aves | | 19.5 | 14 | 8 | 5 | 3 | 12 | 9.5 | 19.5 | 4 | 2 |
| FM of Bests | | 0.085317 | 0.055952 | 0.054365 | 0.025 | 0.035317 | 0.032937 | 0.034921 | 0.09127 | 0.031349 | 0.0163 |
| FMR of Bests | | 19 | 13 | 12 | 3 | 7 | 5 | 6 | 20 | 4 | 1 |

Continue the above Table 6:

| Fun. | Index | GWO | WOA | GTO | DA | SSA | ECBO | MVO | WEO | MKE | LJA |
|---|---|---|---|---|---|---|---|---|---|---|---|
| F24 | Ave | 4.75E00 | 1.36E-03 | 3.1E00 | 1.14E02 | 4.33E01 | 4.67E01 | 7.33E01 | 1.8E-07 | 1.05E02 | 7.18E01 |
| | Std | 5.36E00 | 7.33E-04 | 1.17E00 | 8.57E01 | 5.68E01 | 5.07E01 | 7.85E01 | 3.89E-07 | 9.72E01 | 9.05E01 |
| | Best | 2.36E-05 | 2.53E-04 | 1.39E00 | 4.06E00 | 7.88E-12 | 5.9E-26 | 7.74E-06 | 3.28E-11 | 8.19E-19 | 2.66E-18 |
| F25 | Ave | 1.57E01 | 1.06E01 | 9.25E01 | 1.54E02 | 4.71E01 | 1.19E01 | 5.83E01 | 2.98E00 | 1.28E02 | 1.17E02 |
| | Std | 2.07E01 | 8.67E00 | 4.25E01 | 7.9E01 | 3.73E01 | 2.9E01 | 5.38E01 | 2.17E00 | 8.38E01 | 7.82E01 |
| | Best | 9.52E-03 | 1.39E-01 | 5.79E01 | 4.7E01 | 3.15E00 | 3.63E-27 | 3.16E-03 | 4.82E-01 | 3.51E01 | 5.9E00 |
| F26 | Ave | 4.47E02 | 5.43E02 | 2.41E02 | 1.41E03 | 5.11E02 | 3.09E02 | 4.48E02 | 7.24E00 | 6.8E02 | 8.14E02 |
| | Std | 1.58E02 | 2.72E02 | 1.54E02 | 1.91E02 | 3.75E02 | 3.2E02 | 3.37E02 | 5.02E00 | 3.42E02 | 4.01E02 |
| | Best | 8.62E-03 | 5.13E-02 | 8.13E01 | 9.59E02 | 7E-10 | 1.29E-29 | 3.45E-04 | 6.15E-01 | 2.17E02 | 2.84E02 |
| F27 | Ave | 2.21E01 | 6.31E00 | 6.31E00 | 2.32E01 | 2.64E01 | 6.31E00 | 3.06E01 | 6.31E00 | 3.5E01 | 1.5E01 |
| | Std | 4.93E01 | 8.75E-05 | 1.21E-04 | 2.19E01 | 5.33E01 | 6.56E-10 | 5.67E01 | 1.18E-11 | 5.97E01 | 3.31E01 |
| | Best | 6.32E00 | 6.31E00 | 6.31E00 | 6.31E00 | 6.31E00 | 6.31E00 | 6.31E00 | 6.31E00 | 6.31E00 | 6.31E00 |
| F28 | Ave | 7.71E-01 | 4.51E-01 | 1.89E00 | 9.41E00 | 4.04E01 | 4.38E01 | 7.04E01 | 3.65E-01 | 8.56E01 | 3.95E01 |
| | Std | 7.28E-01 | 3.64E-01 | 2.5E-01 | 7.13E00 | 5.69E01 | 5.83E01 | 6.41E01 | 1.08E-01 | 6.39E01 | 6.82E01 |
| | Best | 3.39E-04 | 1.71E-01 | 1.47E00 | 2.29E00 | 2.81E-01 | 0E00 | 1.34E-01 | 1.7E-01 | 5.88E-01 | 1.02E00 |
| F29 | Ave | 8.05E02 | 7.75E02 | 5.37E02 | 8.42E02 | 4.36E02 | 5.7E02 | 5.47E02 | 2.67E01 | 7.24E02 | 7.51E02 |
| | Std | 3.45E02 | 1.46E02 | 6.65E01 | 5.81E01 | 8.12E01 | 3.42E02 | 3.71E02 | 1.01E02 | 1.41E02 | 2.16E02 |
| | Best | 8.02E02 | 1.12E00 | 4.26E02 | 6.5E02 | 4E02 | 0E00 | 1.49E-02 | 1.73E-06 | 4E02 | 1.38E-06 |
| FM of Aves | | 0.047619 | 0.036111 | 0.038492 | 0.078571 | 0.053968 | 0.044841 | 0.066667 | 0.010714 | 0.080952 | 0.069444 |

| | | | | | | | | | |
|---|---|---|---|---|---|---|---|---|---|
| FMR of Aves | 11 | 6 | 7 | 17 | 13 | 9.5 | 15 | 1 | 18 | 16 |
| FM of Bests | 0.060317 | 0.051984 | 0.073413 | 0.078175 | 0.049603 | 0.020238 | 0.04246 | 0.047222 | 0.059127 | 0.056746 |
| FMR of Bests | 16 | 11 | 17 | 18 | 10 | 2 | 8 | 9 | 15 | 14 |

Table 7 The final rank for the 29 test problems.

| Algorithms | CMA-ES | OTWO | TLBO | LCBO | GMBA | EO | MFO | SO | VLA | AVLA |
|---|---|---|---|---|---|---|---|---|---|---|
| FM of Aves | 0.092118 | 0.058128 | 0.035632 | 0.038424 | 0.056158 | 0.039737 | 0.053941 | 0.058539 | **0.028407** | **0.021839** |
| FMR of Aves | 20 | 15 | 5 | 6 | 13 | 7 | 11 | 16 | **2** | **1** |
| FM of Bests | 0.063136 | 0.055829 | 0.041051 | **0.030213** | 0.050411 | 0.038998 | 0.04376 | 0.06133 | 0.03867 | **0.026437** |
| FMR of Bests | 18 | 15 | 6 | **2** | 11 | 4 | 7 | 17 | 3 | **1** |

Continue the above Table 7:

| Algorithms | GWO | WOA | GTO | DA | SSA | ECBO | MVO | WEO | MKE | LJA |
|---|---|---|---|---|---|---|---|---|---|---|
| FM of Aves | 0.049754 | 0.043678 | 0.032677 | 0.083415 | 0.056897 | 0.0422 | 0.054762 | 0.029228 | 0.065353 | 0.059113 |
| FMR of Aves | 10 | 9 | 4 | 19 | 14 | 8 | 12 | 3 | 18 | 17 |
| FM of Bests | 0.052709 | 0.049425 | 0.050082 | 0.080706 | 0.055419 | 0.039409 | 0.056076 | 0.050082 | 0.052956 | 0.0633 |
| FMR of Bests | 12 | 8 | 9.5 | 20 | 14 | 5 | 16 | 9.5 | 13 | 19 |

## 5.3. Test on five SCNs

This subsection focuses on five distinct types of nodes within supply chain networks: raw material suppliers, manufacturers, wholesalers, retailers, and markets. Each supplier provides only one type of raw material, while each manufacturer produces a single, standardized product. The first four supply chain networks consist of four levels, while the fifth network extends to five levels. Before presenting the results obtained from 20 algorithms, it is essential to first define and explain the specific implementations of the abstract functions introduced in Section 3 within the context of the Supply Chain Network Evaluation Model (SCNEM).

Let $r_{i,k}$ be the fixed rate of material $k$ to produce $1 + r_{i,T}$ unit of the standardized product at manufacturer $i$ where $r_{i,T} > -1$ is a given constant called transforming rate. The relationship $\sum_k r_{i,k} = 1$ holds. Let $Q_i^P$ be the actually produced product at manufacturer $i$. The manufacture function $[Q_i^O, Q_i^H] = M(Q_i^E)$, $\forall i \in I \setminus I_M$ is realized with the following concrete functions:

$$Q_i^P = (1 + r_{i,T}) \min_k \left\{ \frac{Q_{i,k}^E}{r_{i,k}} \middle| k \in PK_S \right\} \tag{33}$$

$$Q_i^O + Q_i^H = Q_i^P \tag{34}$$

$$Q_{i,k}^H = Q_{i,k}^E - \frac{Q_i^P}{(1+r_{i,T})} r_{i,k}, k \in PK_S \tag{35}$$

At a market $i \in I_M$, the pricing function $p_{i,k}^E = P(Q_{i,k}^E), \forall i \in I_M$ has the following concrete form:

$$P(Q_i^E) = \max \left\{ 0, p_i^{Max} - a_i Q_i^E - b_i Q_i^{E^2} \right\} \tag{36}$$

In equation (36), $p_i^{Max}$ is the maximal price of the standardized product at market $i$ and $a_i$ and $b_i$ are two given parameters. The subscript k is omitted in equation (36) because there is only one standardized product to be considered. The transportation cost associated with flow $f_{i,j,k}$ on link $(i,j)$ is calculated as follows:

$$C_{i,j,k}(f_{i,j,k}) = a_{i,j,k} f_{i,j,k} + b_{i,j,k} f_{i,j,k}^2 + c_{i,j,k} \tag{37}$$

In function (37), $a_{i,j,k}$, $b_{i,j,k}$ and $c_{i,j,k}$ are three given parameters. The general production cost $C_{i,k}(Q_i)$ consists of purchase expenditure, fixed cost, variable production cost, holding cost and transaction cost. The fixed cost is given by a positive constant. We assume that the variable production cost, holding cost and transaction cost all have the following common form:

$$c(Q) = aQ + bQ^2 \tag{38}$$

In function (38), $Q$ is the quantity of the product in question and $a$ and $b$ are two given parameters.

This subsection examines five supply chain networks. The first network, denoted as SCN1, consists of 2 material suppliers, 1 manufacturer, 2 retailers, and 2 markets. The second network, SCN2, has the same structure as SCN1 but includes a key difference: manufacturers in SCN2 can directly supply products to markets. The third network, SCN3, includes 4 material suppliers, 2 manufacturers, 3 retailers, and 4 markets. The fourth network, SCN4, comprises 4 material suppliers, 2 manufacturers, 4 retailers, and 2 markets. Finally, the fifth network, SCN5, extends SCN4 by adding a wholesaler as an additional decision-making entity. The number of variables involved in each network is as follows: SCN1 has 15 variables, SCN2 has 17, SCN3 has 39, SCN4 has 38, and SCN5 has 41.

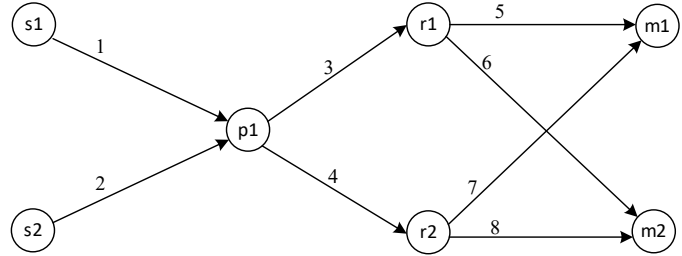

Fig. 1. A multitiered supply chain with one central manufacturer.

The topological structure of SCN1 is illustrated in Figure 1. In the figure, the two suppliers are represented as s1 and s2, the sole manufacturer as p1, the two retailers as r1 and r2, and the two markets as m1 and m2. In SCN1, s1 supplies one type of raw material, while s2 provides another. The upper bounds of variables include $Q_{i,k}^{E,Max} = 500, \forall i \in I_S, k$, $\lambda_{i,k}^{Max} = 1.0, \forall i \in I \backslash I_M, k$ and $f_{i,j,k}^{Max} = 5000, \forall i,j,k$ in SCN1.

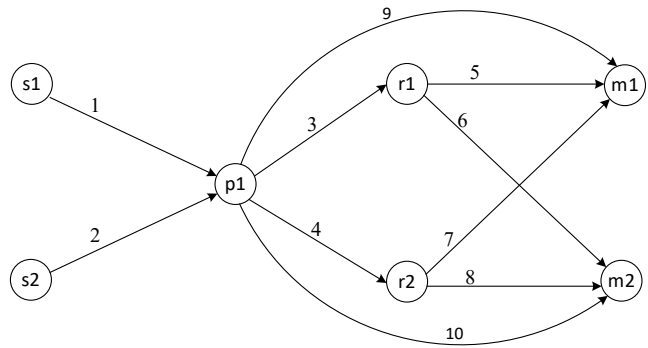

Fig. 2. A multitiered supply chain with direct transactions between manufacturer and markets.

The topologic structure of SVN2 is given in Figure2. Comparing to SCN1, SCN2 adds two connectors from the manufacturer p1 to two markets m1 and m2. The other things are same between SCN1 and SCN2.

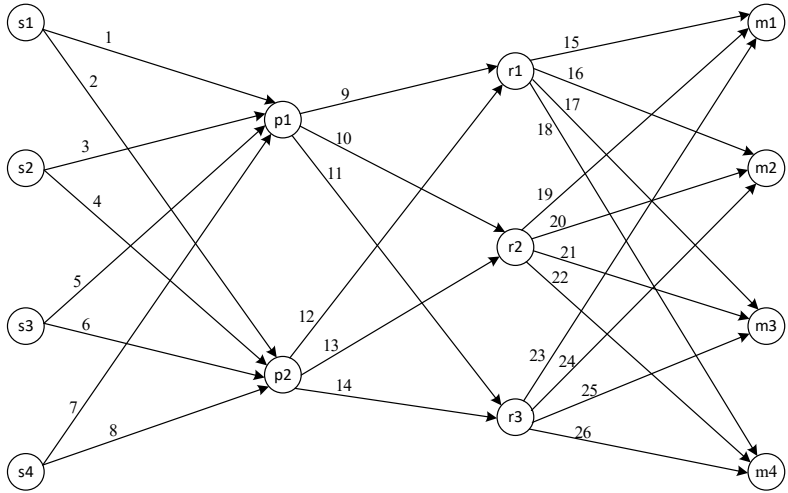

Fig. 3. The topological network of the third supply chain network.

SCN3 is more complicated than SCN1 and SCN2 in terms of the number of decision makers and the number of constraints. Figure 3 depicts the topological structure of SCN3. In Figure 3, four raw material suppliers (labeled s1, s2, s3, and s4 in Fig. 2) provide two types of raw materials; raw material one is supplied by s1 and s2, and raw material two by s3 and s4. Two manufacturers are represented as p1 and p2 in Fig. 3. Three retailers (indicated by r1, r2, and r3) supply finished products to four markets (indicated by m1, m2, m3, and m4). The related upper bounds on variables are given by $Q_{i,k}^{E,Max} = 10, \forall i \in I_S, k$, $\lambda_{i,k}^{Max} = 1.0, \forall i \in I \backslash I_M, k$ and $f_{i,j,k}^{Max} = 5000, \forall i,j,k$ for SCN2.

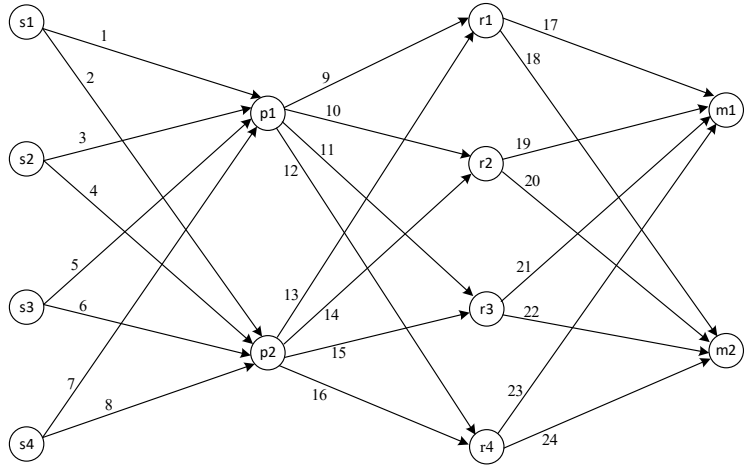

Fig. 4. The topological network of the fourth supply chain. network.

SCN4 is depicted in Figure 4. The four levels of SCN4 correspond to material suppliers, manufacturers, retailers and markets, respectively. There are total 24 links in SCN4. Four suppliers of raw material are labeled with s1, s2, s3, and s4, respectively. In SCN4, s1 and s2 supply one raw material and s3 and s4 the other raw material. Two manufacturers are indicated by p1 and p2 in Figure 4. Four retailers are labeled with r1, r2, r3, and r4, respectively. Two markets m1 and m2 are at the last level of SCN4. The related upper bounds to variables are given by $Q_{i,k}^{E,Max} = 50, \forall i \in I_S, k$, $\lambda_{i,k}^{Max} = 1.0, \forall i \in I \backslash I_M, k$ and $f_{i,j,k}^{Max} = 5000, \forall i,j,k$.

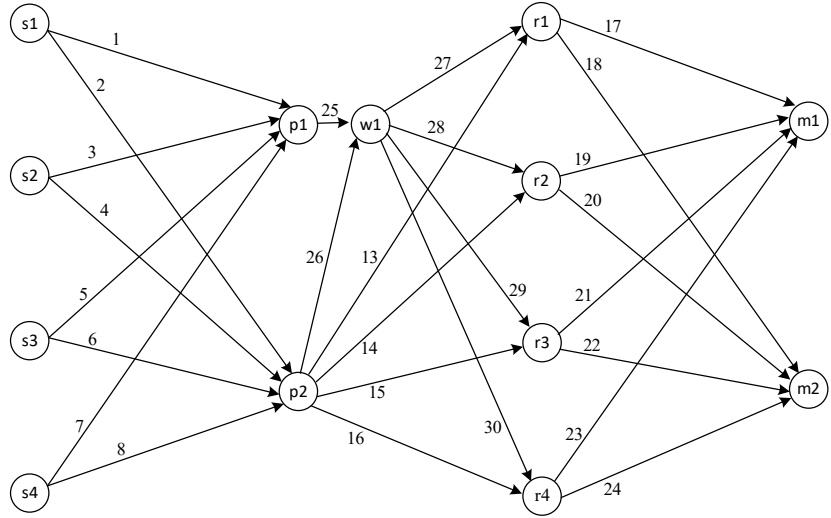

Fig. 5. The topological network of the fifth supply chain. network.

The fifth supply chain network (SCN5), depicted in Figure 5, is a five-level network. Unlike SCN4, SCN5 introduces a wholesaler, denoted as w1. All products produced by p1 are shipped to w1, which also receives products from p2. Both w1 and p2 have the capability to sell products directly to retailers. In SCN5, links 9 through 12 from SCN4 are removed, and new links numbered 25 through 30 are added. It is assumed that any links shared between SCN4 and SCN5 with the same series number have identical transportation cost functions.

Appendix B provides detailed initialization information for these supply chain networks. Additionally, five solutions generated by AVLA for SCN1 through SCN5 are included in the appendix. This information allows interested readers to replicate the experiments and verify the solutions.

The primary results obtained from 20 heuristic algorithms across the five SCNs are summarized in Table 8. The data indicate that AVLA outperforms the other 19 algorithms based on the mean, standard deviation, and best objective values. AVLA consistently achieves the top ranks in both FMR of Aves and FMR of Bests. Since the theoretical optimal objective values for these SCNs are zero, solutions closer to zero are considered superior. A detailed comparison of the algorithms reveals that AVLA produces significantly better-quality solutions than the others. VLA and WEO secure the second and third ranks in the FMR of Aves for three SCNs, while TLBO and LCBO occupy the second and third ranks in the FMR of Bests, respectively. Although the rankings of algorithms in these SCNs differ slightly from those presented in Table 7 for 29 test benchmarks, the consistent and outstanding performance of AVLA remains unchanged.

**Table 8** Optimization results and comparison for composite benchmark functions.

| Fun. | Index | CMA-ES | OTWO | TLBO | LCBO | GMBA | EO | MFO | SO | VLA | AVLA |
|---|---|---|---|---|---|---|---|---|---|---|---|
| SCN1 | Ave | 7E05 | 4.83E01 | 3.53E-01 | 9.52E-01 | 6.8E01 | 2E00 | 9.47E00 | 9.69E01 | 7.75E-02 | 0 |
|  | Std | 3.54E05 | 3.75E01 | 1.65E00 | 3.75E00 | 8.2E01 | 3.17E00 | 1.18E01 | 2.3E01 | 6.37E-02 | 0 |
|  | Best | 2.75E05 | 1.51E01 | 0 | 0 | 8.22E00 | 0 | 0 | 2.11E01 | 4.9E-03 | 0 |
| SCN2 | Ave | 1.75E06 | 1.22E03 | 1.46E02 | 1.32E03 | 8.72E02 | 1.4E02 | 5.44E02 | 9.5E01 | 1.94E02 | 5.38E-01 |
|  | Std | 6.74E05 | 2.51E03 | 1.5E02 | 6.04E03 | 1.04E03 | 2.97E02 | 8.11E02 | 2.02E01 | 5.86E01 | 1.25E00 |
|  | Best | 9.16E04 | 2.16E01 | 3.13E00 | 1.16E-04 | 4.19E01 | 1.04E01 | 5.31E00 | 4.03E01 | 6.36E01 | 0 |
| SCN3 | Ave | 4.86E03 | 6.17E02 | 2.17E02 | 3.12E02 | 2.83E02 | 1.9E02 | 3.89E02 | 1.11E03 | 9.21E01 | 1.25E00 |
|  | Std | 5.36E03 | 2.48E02 | 1.16E02 | 2.5E02 | 1.22E02 | 1.62E02 | 2.38E02 | 1.81E02 | 2.35E01 | 1.01E00 |
|  | Best | 5.11E01 | 2.27E02 | 1.12E01 | 5.67E00 | 3.57E01 | 2.77E01 | 1.68E01 | 7.08E02 | 4.79E01 | 2.22E-01 |
| SCN4 | Ave | 4.53E05 | 1.13E04 | 8.79E02 | 1.9E03 | 8.29E03 | 1.01E03 | 2.17E03 | 1.11E04 | 6.52E02 | 2.36E01 |
|  | Std | 2.08E05 | 1.1E04 | 6.21E02 | 1.64E03 | 5.33E03 | 6.74E02 | 1.2E03 | 2.47E03 | 2.16E02 | 5.92E01 |
|  | Best | 3.94E04 | 1.67E03 | 2.25E01 | 1.74E02 | 3.14E03 | 6.21E01 | 4.15E02 | 7.12E03 | 3.98E02 | 1.96E-02 |
| SCN5 | Ave | 7.46E04 | 1.43E03 | 5.79E01 | 1.09E02 | 1.07E03 | 7.62E01 | 8.81E01 | 2.37E02 | 2.21E01 | 6.58E00 |
|  | Std | 7E04 | 4.06E03 | 2.76E01 | 5.23E01 | 2.8E03 | 4.43E01 | 6.25E01 | 1.66E01 | 6.41E00 | 5.86E00 |
|  | Best | 1.78E03 | 8.73E01 | 8.74E00 | 3.84E01 | 6.9E01 | 1.76E01 | 1.51E00 | 1.82E02 | 8.17E00 | 5.39E-01 |

|  | FM of Aves | 0.088571 | 0.069524 | 0.024762 | 0.045714 | 0.059048 | 0.026667 | 0.046667 | 0.059048 | **0.019048** | **0.004762** |
|---|---|---|---|---|---|---|---|---|---|---|---|
|  | FMR of Aves | 19 | 16.5 | 4 | 7 | 13.5 | 5.5 | 8.5 | 13.5 | **2** | **1** |
|  | FM of Bests | 0.08381 | 0.063333 | **0.020476** | 0.023333 | 0.06 | 0.032857 | 0.024286 | 0.077143 | 0.045714 | **0.007143** |
|  | FMR of Bests | 19 | 14 | **2** | 3 | 12 | 7 | 4 | 17 | 10 | **1** |

Continue the above Table 8:

| Fun. | Index | GWO | WOA | GTO | DA | SSA | ECBO | MVO | WEO | MKE | LJA |
|---|---|---|---|---|---|---|---|---|---|---|---|
| SCN1 | Ave | 1.32E06 | 1.95E03 | 7.02E00 | 2.89E02 | 6.24E03 | 1.7E05 | 9.56E01 | 2.68E00 | 8.95E00 | 2.07E03 |
|  | Std | 3.34E05 | 1.06E04 | 2.23E00 | 1.79E02 | 2.27E04 | 6.14E04 | 1.45E02 | 1.76E00 | 1.17E01 | 7.66E03 |
|  | Best | 8.02E05 | 2.12E-02 | 3.84E00 | 6.49E01 | 1.87E01 | 1.29E04 | 1.09E01 | 7.55E-01 | 0 | 1.7E-09 |
| SCN2 | Ave | 2.13E06 | 2.34E01 | 7.35E02 | 2.87E02 | 4.54E03 | 3.24E01 | 2.61E02 | 6.45E01 | 4.61E02 | 4.06E02 |
|  | Std | 2.93E05 | 4.76E00 | 2.16E03 | 1.51E02 | 3.56E03 | 5.1E01 | 3.24E02 | 9.44E01 | 8.42E02 | 4.95E02 |
|  | Best | 1.55E06 | 2.03E01 | 2.39E01 | 6.32E01 | 8.21E01 | 5.5E-01 | 4.4E01 | 7.02E00 | 1E01 | 1.62E-01 |
| SCN3 | Ave | 3.17E06 | 5.19E05 | 3.95E02 | 1.08E03 | 3E02 | 4.8E01 | 7.78E01 | 7.47E01 | 4.46E02 | 4.67E02 |
|  | Std | 8.04E04 | 8.4E05 | 6.7E01 | 2.34E02 | 1.06E02 | 2.77E01 | 6.9E01 | 2.57E01 | 2.64E02 | 1.56E02 |
|  | Best | 2.8E06 | 4.5E02 | 2.42E02 | 6.42E02 | 7.89E01 | 6.23E00 | 2.09E01 | 2.49E01 | 2.07E01 | 2.19E02 |
| SCN4 | Ave | 8.36E05 | 3.1E05 | 2.62E03 | 1.5E04 | 2.12E04 | 6.2E02 | 3.47E03 | 1.87E03 | 2.87E03 | 1.98E03 |
|  | Std | 0E00 | 1.36E05 | 5.74E02 | 3.79E03 | 1.71E04 | 4.73E02 | 1.56E03 | 6.53E02 | 2.18E03 | 9.82E02 |
|  | Best | 8.36E05 | 4.34E03 | 1.67E03 | 7.28E03 | 4.12E03 | 1.8E00 | 1.43E03 | 9.16E02 | 1.03E01 | 1.46E02 |
| SCN5 | Ave | 8.58E05 | 2.36E05 | 9.22E01 | 1.37E03 | 5.93E02 | 4.84E01 | 2.26E02 | 5.63E01 | 9.6E01 | 1.27E02 |
|  | Std | 7.82E04 | 2.19E05 | 1.35E01 | 5.98E02 | 7.77E02 | 5.94E02 | 1.52E01 | 1.52E01 | 7.95E01 | 2.04E02 |
|  | Best | 6.92E05 | 2.34E02 | 5.95E01 | 4.12E02 | 9.45E01 | 4.49E00 | 4.56E01 | 2.67E01 | 1.47E01 | 1.15E01 |
| FM of Aves |  | 0.095238 | 0.069524 | 0.048571 | 0.068571 | 0.071429 | 0.026667 | 0.046667 | 0.021905 | 0.050476 | 0.057143 |
| FMR of Aves |  | 20 | 16.5 | 10 | 15 | 18 | 5.5 | 8.5 | 3 | 11 | 12 |
| FM of Bests |  | 0.095238 | 0.065714 | 0.060476 | 0.082857 | 0.072381 | 0.028571 | 0.054286 | 0.041905 | 0.02619 | 0.034286 |
| FMR of Bests |  | 20 | 15 | 13 | 18 | 16 | 6 | 11 | 9 | 5 | 8 |

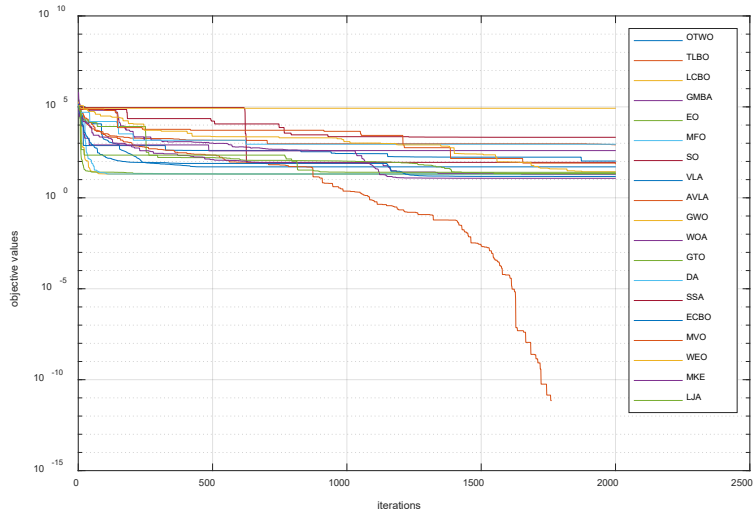

Fig. 6 The trends of objective values of SCN2 over iterations for 19 algorithms.

Figure 6 illustrates the changes in objective values for SCN2 over increasing iterations, as obtained by 19 algorithms. For each algorithm, the data represent a random run used to solve the equilibrium model of SCN2. CMA-ES is excluded from Figure 6 because it consistently terminates prematurely before reaching the maximum number of iterations and yields relatively large objective values for SCN2. The final best objective function values for the 19 algorithms—OTWO, TLBO, LCBO, GMBA, EO, MFO, SO, VLA, AVLA, GWO, WOA, GTO, DA, SSA, ECBO, MVO, WEO, MKE, and LJA—are 50.3808, 820.8218, 20.2251, 394.5163, 20.4178, 20.3930, 88.8839, 105.2374, 0, 83881.1589, 20.3135, 21.6553, 881.6225, 2125.4000, 14.9359, 73.9761, 27.0899, 11.8126, and 25.2408, respectively. Among these, AVLA achieves the theoretical optimum of 0 at iteration 1762. The curves in Figure 6 reveal that AVLA surpasses all other algorithms, providing the best solution for SCN2 starting from iteration 900 onward.

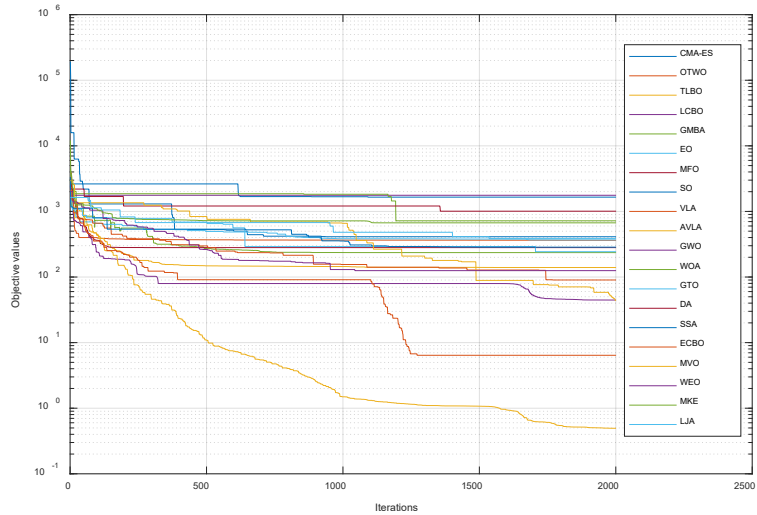

Fig. 7. The trends of objective values of SCN3 over iterations for 20 algorithms.

Figure 7 illustrates the trends of objective function values for SCN3 over iterations across all 20 algorithms. For each algorithm, the data are derived from a random run to solve the equilibrium model of SCN3. While the initial values are closely grouped for most algorithms, the final values vary significantly over a wider range. AVLA achieves the smallest final objective value of 0.49395, outperforming all other algorithms. The second-best final value, 6.41012, is obtained by ECBO, followed by MVO and LCBO with final values of 44.12364 and 44.52552, respectively. In contrast, GWO produces the worst final objective value of 1763.71717. The curves in Figure 7 reveal that AVLA consistently yields lower objective values than the other 19 algorithms after just 200 iterations. This result highlights AVLA's capability to achieve superior solutions even with a relatively small number of iterations.

To conclude this section, we clarify why a fixed population size and a fixed maximum number of iterations were chosen as the predefined parameters for comparing algorithms.

In the literature, alternative comparison conditions, such as fixed computational time or fixed Number of Function Evaluations (NFE), are also commonly used. Table 9 summarizes the mean computational times and mean NFEs for each algorithm applied to SCN4 over 30 runs under our fixed-parameter setting. The data reveal that most algorithms require approximately 100,000 objective function evaluations per run, with variations arising from their internal mechanisms. For instance, AVLA performs 112,001 evaluations due to additional operations like reflection-based estimations. Similarly, algorithms often modify randomly generated solutions to ensure feasibility, resulting in increased NFEs.

Using fixed NFEs as a comparison metric overlooks other computational operations, such as the sorting operations inherent in many algorithms. Table 9 also shows that even with similar NFEs, computational times can vary significantly. Notably, a higher NFE does not necessarily correspond to longer computational times. These findings demonstrate that neither fixed NFEs nor fixed computational times are suitable comparison metrics for our study.

Another reason for selecting a fixed population size of 50 and a maximum of 2000 iterations is to allow each algorithm to fully demonstrate its potential in solving complex problems. The smaller NFE observed for CMA-ES arises from its mechanism, which forces termination when further improvement is no longer possible.

Table 9 Average computational Time and average NFEs(millisecond-ms) of SCN4.

| Algorithms | CMA-ES | OTWO | TLBO | LCBO | GMBA | EO | MFO | SO | VLA | AVLA |
|---|---|---|---|---|---|---|---|---|---|---|
| Time(ms) | 1037 | 10132 | 4222 | 3988 | 2326 | 3012 | 3037 | 2507 | 2766 | 3970 |
| NFEs | 18966 | 143702 | 200050 | 198662 | 100000 | 102051 | 100050 | 101696 | 120324 | 112001 |

Continue the above Table 9:

| Algorithms | GWO | WOA | GTO | DA | SSA | ECBO | MVO | WEO | MKE | LJA |
|---|---|---|---|---|---|---|---|---|---|---|
| Time(ms) | 2984 | 2953 | 6136 | 5051 | 2855 | 3029 | 3302 | 2187 | 3428 | 6061 |
| NFEs | 100050 | 100050 | 500050 | 178855 | 98050 | 109976 | 100050 | 100050 | 200050 | 100050 |

## 6. Conclusions

In this paper, we propose a generalized Supply Chain Network Equilibrium Model (SCNEM) capable of addressing a wide range of supply chain configurations. By employing abstract production and demand market pricing functions, the SCNEM can accommodate complex real-world scenarios. These abstract functions can be replaced with concrete ones when applied to specific supply chains.

To efficiently solve the highly nonlinear, non-convex, and non-smooth SCNEM, we introduce a novel heuristic algorithm, AVLA. Inspired by group learning behaviors, AVLA simulates the learning processes of individuals within different subgroups at various stages. It leverages a success history-based parameter adaptation method to reduce the manual tuning effort. Extensive computational experiments on 29 benchmark problems and 5 SCNs demonstrate AVLA's superior performance compared to 19 state-of-the-art algorithms, often providing the only feasible solutions for complex supply chain equilibrium problems.

Future research directions include extending the SCNEM to specialized supply chains such as emergency material and fresh-food supply chains, and incorporating environmental factors like pollution. While AVLA's adaptive mechanism reduces parameter tuning, further refinement is needed to optimize memory usage and proper limits to non-improvement iterations. Additionally, investigating the conditions under which adaptive AVLA outperforms its non-adaptive counterpart will be a valuable research avenue.

**Disclosure statement**
No potential conflict of interest was reported by the authors.

**Data availability statement**
The data that support the findings of this study are available from the corresponding author upon request.

# Appendix A (Details of the 29 optimization problems):

**Table 1** Unimodal benchmark functions.

| Functions | Dim | Range | $f_{min}$ |
|---|---|---|---|
| $f_1(x) = \sum_{i=1}^{n} x_i^2$ | 30 | $[-100,100]^n$ | $f(0,0,\dots,0) = 0$ |
| $f_2(x) = \sum_{i=1}^{n}|x_i| + \prod_{i=1}^{n}|x_i|$ | 30 | $[-10,10]^n$ | $f(0,0,\dots,0) = 0$ |
| $f_3(x) = \sum_{i=1}^{n}(\sum_{j=1}^{i} x_j)^2$ | 30 | $[-100,100]^n$ | $f(0,0,\dots,0) = 0$ |
| $f_4(x) = \max_{i}\{|x_i|, 1 \leq i \leq n\}$ | 30 | $[-100,100]^n$ | $f(0,0,\dots,0) = 0$ |
| $f_5(x) = \sum_{i=1}^{n-1}[100(x_{i+1} - x_i^2)^2 + (x_i - 1)^2]$ | 30 | $[-30,30]^n$ | $f(1,1,\dots,1) = 0$ |
| $f_6(x) = \sum_{i=1}^{n}([x_i + 0.5])^2$ | 30 | $[-100,100]^n$ | 0 |
| $f_7(x) = \sum_{i=1}^{n} i x_i^4 + random[0,1)$ | 30 | $[-1.28,1.28]^n$ | 0 |

**Table 2** Multimodal benchmark functions.

| Functions | Dim | Range | $f_{min}$ |
|---|---|---|---|
| $f_8(x) = \sum_{i=1}^{n}\left[-x_i \sin\left(\sqrt{|x_i|}\right)\right]$ | 30 | $[-500,500]^n$ | $-418.9829n$ |
| $f_9(x) = \sum_{i=1}^{n}[x_i^2 - 10\cos(2\pi x_i) + 10]$ | 30 | $[-5.12,5.12]^n$ | 0 |
| $f_{10}(x) = -20\exp\left(-0.2\sqrt{\frac{1}{n}\sum_{i=1}^{n} x_i^2}\right) - \exp(\frac{1}{n}\sum_{i=1}^{n}\cos(2\pi x_i)) + 20 + e$ | 30 | $[-32,32]^n$ | 0 |
| $f_{11}(x) = \sum_{i=1}^{d}\frac{x_i^2}{4000} - \prod_{i=1}^{d}\cos\left(\frac{x_i}{\sqrt{i}}\right) + 1$ | 30 | $[-512,512]^n$ | 0 |
| $f_{12}(x) = \frac{\pi}{n}\{10\sin^2(\pi y_1) + \sum_{i=1}^{n-1}(y_i - 1)^2[1 + 10\sin^2(\pi y_{i+1})] + (y_n - 1)^2\} + \sum_{i=1}^{n} u(x_i,10,100,4)$ | 30 | $[-50,50]^n$ | 0 |
| $y_i = 1 + \frac{x_i+1}{4},\ u(x_i, a, k, m) = \begin{cases} k(x_i - a)^m & x_i > a \\ 0 & -a \leq x_i \leq a \\ k(-x_i - a)^m & x_i < -a \end{cases}$ | | | |
| $f_{13}(x) = 0.1\left\{\sin^2(3\pi x_1) + \sum_{i=1}^{n}(x_i - 1)^2[1 + \sin^2(3\pi x_i + 1)] + (x_n - 1)^2[1 + \sin^2(2\pi x_n)]\right\} + \sum_{i=1}^{n} u(x_i,5,100,4)$ | 30 | $[-50,50]^n$ | 0 |

**Table 3** Fixed dimension multimodal benchmark functions.

| Functions | Dim | Range | $f_{min}$ |
|---|---|---|---|
| $f_{14}(x) = \left(\frac{1}{500} + \sum_{j=1}^{25}\frac{1}{j+\sum_{i=1}^{2}(x_i - a_{ij})^6}\right)^{-1}$ | 2 | $[-65.536,65.536]^n$ | 1 |
| $f_{15}(x) = \sum_{i=1}^{11}\left[a_i - \frac{x_1(b_i^2 + b_i x_2)}{b_i^2 + b_i x_3 + x_4}\right]^2$ | 4 | $[-5,5]^n$ | 0.0003075 |
| $f_{16}(x) = 4x_1^2 - 2.1x_1^4 + x_1^6/3 + x_1 x_2 - 4x_2^2 + 4x_2^4$ | 2 | $[-5,5]^n$ | -1.0316285 |
| $f_{17}(x) = \left(x_2 - \frac{5.1}{4\pi^2}x_1^2 + \frac{5}{\pi}x_1 - 6\right)^2 + 10\left(1 - \frac{1}{8\pi}\right)\cos x_1 + 10$ | 2 | $[-5,5]^n$ | 0.398 |
| $f_{18}(x) = [1 + (x_1 + x_2 + 1)^2(19 - 14x_1 + 3x_1^2 - 14x_1 + 6x_1 x_2 + 3x_2^2)] \times [30 + (2x_1 - 3x_2)^2(18 - 32x_1 + 12x_1^2 + 48x_2 - 36x_1 x_2 + 27x_2^2)]$ | 2 | $[-2,2]^n$ | 3 |
| $f_{19}(x) = -\sum_{i=1}^{4} c_i \exp\left(-\sum_{j=1}^{3} a_{ij}(x_j - p_{ij})^2\right)$ | 3 | $[0,1]^n$ | -3.86 |
| $f_{20}(x) = -\sum_{i=1}^{4} c_i \exp\left(-\sum_{j=1}^{6} a_{ij}(x_j - p_{ij})^2\right)$ | 6 | $[0,1]^n$ | -3.32 |
| $f_{21}(x) = -\sum_{i=1}^{5}[(x - a_i)(x - a_i)^T + c_i]^{-1}$ | 4 | $[0,10]^n$ | -10.1532 |
| $f_{22}(x) = -\sum_{i=1}^{7}[(x - a_i)(x - a_i)^T + c_i]^{-1}$ | 4 | $[0,10]^n$ | -10.4028 |
| $f_{23}(x) = -\sum_{i=1}^{10}[(x - a_i)(x - a_i)^T + c_i]^{-1}$ | 4 | $[0,10]^n$ | -10.5363 |

**Table 4** Composite benchmark functions.

| Functions | Dim | Range | $f_{min}$ |
|---|---|---|---|
| $F_{24}$ (CF1) <br> $f_1, f_2, f_3, \ldots, f_{10}$ =Sphere Function <br> $[\sigma_1, \sigma_2, \sigma_3, \ldots, \sigma_{10}] = [1,1,1,\ldots,1]$ <br> $[\lambda_1, \lambda_2, \lambda_3, \ldots, \lambda_{10}] = [5/100, 5/100, 5/100, \ldots, 5/100]$ | 10 | [-5,5] | 0 |
| $F_{25}$ (CF2) <br> $f_1, f_2, f_3, \ldots, f_{10}$ =Griewank's Function <br> $[\sigma_1, \sigma_2, \sigma_3, \ldots, \sigma_{10}] = [1,1,1,\ldots,1]$ <br> $[\lambda_1, \lambda_2, \lambda_3, \ldots, \lambda_{10}] = [5/100, 5/100, 5/100, \ldots, 5/100]$ | 10 | [-5,5] | 0 |
| $F_{26}$ (CF3) <br> $f_1, f_2, f_3, \ldots, f_{10}$ =Griewank's Function <br> $[\sigma_1, \sigma_2, \sigma_3, \ldots, \sigma_{10}] = [1,1,1,\ldots,1]$ <br> $[\lambda_1, \lambda_2, \lambda_3, \ldots, \lambda_{10}] = [1,1,1,\ldots,1]$ | 10 | [-5,5] | 0 |
| $F_{27}$ (CF4) <br> $f_1, f_2$= Ackley's Function, $f_3, f_4$= Rastrigin's Function, <br> $f_5, f_6$= Weierstrass Function, $f_7, f_8$= Griewank's Function, <br> $f_9, f_{10}$= Sphere Function <br> $[\sigma_1, \sigma_2, \sigma_3, \ldots, \sigma_{10}] = [1,1,1,\ldots,1]$ <br> $[\lambda_1, \lambda_2, \lambda_3, \ldots, \lambda_{10}]$ =[5/32, 5/32, 1, 1, 5/0.5, 5/0.5, 5/100, 5/100, 5/100, 5/100] | 10 | [-5,5] | 0 |
| $F_{28}$ (CF5) <br> $f_1, f_2$= Rastrigin s Function, $f_3, f_4$= Weierstrass Function, <br> $f_5, f_6$= Griewank s Function, $f_7, f_8$= Ackley s Function, <br> $f_9, f_{10}$ = Sphere Function <br> $[\sigma_1, \sigma_2, \sigma_3, \ldots, \sigma_{10}] = [1,1,1,\ldots,1]$ <br> $[\lambda_1, \lambda_2, \lambda_3, \ldots, \lambda_{10}]$ =[1/5, 1/5, 5/0.5, 5/0.5, 5/100, 5/100, 5/32, 5/32, 5/100, 5/100] | 10 | [-5,5] | 0 |
| $f_{29}$ (CF6) <br> $f_1, f_2$= Rastrigin's Function, $f_3, f_4$= Weierstrass' Function, <br> $f_5, f_6$= Griewank's Function, $f_7, f_8$= Ackley's Function, <br> $f_9, f_{10}$= Sphere Function <br> $[\sigma_1, \sigma_2, \sigma_3, \ldots, \sigma_{10}]$ =[0.1, 0.2, 0.3, 0.4, 0.5, 0.6, 0.7, 0.8, 0.9, 1] <br> $[\lambda_1, \lambda_2, \lambda_3, \ldots, \lambda_{10}]$ =[0.1 ∗ 1/5, 0.2 ∗ 1/5, 0.3 ∗ 5/0.5, 0.4 ∗ 5/0.5, 0.5 ∗ 5/100, <br> 0.6 ∗ 5/100, 0.7 ∗ 5/32, 0.8 ∗ 5/32, 0.9 ∗ 5/100, 1 ∗ 5/100] | 10 | [-5,5] | 0 |

# Appendix B (Information and solutions to five SCNs):

*Initialization data of SCN 1:*

**Table 1** Parameters related to costs for SCN1.

| Spots | FC | VC | | HC | | TC | |
|---|---|---|---|---|---|---|---|
| | | $a$ | $b$ | $a$ | $b$ | $a$ | $b$ |
| s1 | 34 | 0.02 | 0.000012 | 0.01 | 0.0 | 0.002 | 0.0 |
| s2 | 41 | 0.01 | 0.00001 | 0.02 | 0.0 | 0.002 | 0.0 |
| p1 | 10 | 0.001 | 0.000003 | 0.001 | 0.0 | 0.0001 | 0.0 |
| r1 | 10 | 0.0001 | 0.000001 | 0.005 | 0.0 | 0.005 | 0.0 |
| r2 | 20 | 0.0002 | 0.000001 | 0.004 | 0.0 | 0.0006 | 0.0 |

**Table 2** Parameters related to production process and holding cost for SCN1.

| Manufacturer $i$ | $r_{i,1}$ | $r_{i,2}$ | $r_{i,T}$ | HC of raw material 1 | | HC of raw material 2 | |
|---|---|---|---|---|---|---|---|
| | | | | $a$ | $b$ | $a$ | $b$ |
| p1 | 0.3 | 0.7 | 1.0 | 0.001 | 0.0 | 0.001 | 0.0 |

**Table 3** Parameters associated with pricing function for SCN1.

| Market $m$ | $p_m^{Max}$ | $a_m$ | $b_m$ |
|---|---|---|---|
| m1 | 82.9 | 0.0032 | 0.000076 |
| m2 | 92.8 | 0.004 | 0.000045 |

**Table 4** The parameters of transportation cost function for SCN1.

| No.$(i,j)$ | $a_{ij}$ | $b_{ij}$ | $c_{ij}$ |
|---|---|---|---|
| 1 | 0.0004 | 0.00002 | 0.5 |
| 2 | 0.0004 | 0.00003 | 0.5 |
| 3 | 0.0005 | 0.00002 | 0.5 |
| 4 | 0.0004 | 0.00004 | 0.5 |
| 5 | 0.0004 | 0.00005 | 0.5 |
| 6 | 0.0006 | 0.00002 | 0.5 |
| 7 | 0.0004 | 0.00006 | 0.5 |
| 8 | 0.0008 | 0.00002 | 0.5 |

*One solution to SCN1 with objective value 0 is given below:*

**Table 5** Equilibrium results of SCN1 associated with links.

| No. $(i,j)$ | $p_i^O$ | $c_{ij}$ | $p_j^E$ | $p_j^E - c_{ij} - p_i^O$ | $f_{i,j}$ |
|---|---|---|---|---|---|
| 1 | 93.9405 | 0.5002 | 94.4408 | 0.0000 | 0.4875 |
| 2 | 63.6962 | 0.5006 | 64.1968 | 0.0000 | 1.1375 |
| 3 | 54.7340 | 0.5011 | 55.2350 | 0.0000 | 1.5566 |
| 4 | 54.7340 | 0.5010 | 55.2350 | 0.0000 | 1.6933 |
| 5 | 92.2858 | 0.5000 | 82.9000 | -9.8858 | 0.0000 |
| 6 | 92.2858 | 0.5011 | 92.7869 | 0.0000 | 1.5566 |
| 7 | 92.2853 | 0.5000 | 82.9000 | -9.8853 | 0.0000 |
| 8 | 92.2853 | 0.5016 | 92.7869 | 0.0000 | 1.6933 |

**Table 6** Equilibrium results of SCN1 associated with material suppliers.

| Supplier $s$ | $Q_s^E$ | $Q_s^O$ | $Q_s^H$ | $p_s^O$ | $\lambda_s$ |
|---|---|---|---|---|---|
| s1 | 0.4875 | 0.4875 | 0.0000 | 93.9405 | 0.3465 |
| s2 | 1.1375 | 1.1375 | 0.0000 | 63.6962 | 0.7666 |

**Table 7** Equilibrium results of SCN1 associated with manufacturers.

| Manufacturer $i$ | $Q_i^P$ | $Q_i^O$ | $Q_i^H$ | $Q_{i,s1}^H$ | $Q_{i,s2}^H$ | $p_i^O$ | $p_{i,s1}^E$ | $p_{i,s2}^E$ | $\lambda_i$ |
|---|---|---|---|---|---|---|---|---|---|
| p1 | 3.2499 | 3.2499 | 0.0000 | 0.0000 | 0.0000 | 54.7340 | 94.4408 | 64.1968 | 0.3782 |

**Table 8** Equilibrium results of SCN1 associated with retailers.

| Retailer $j$ | $Q_j^E$ | $Q_j^O$ | $Q_j^H$ | $p_j^E$ | $p_j^O$ | $\lambda_j$ |
|---|---|---|---|---|---|---|
| r1 | 1.5566 | 1.5566 | 0.0000 | 55.2350 | 92.2858 | 0.4966 |
| r2 | 1.6933 | 1.6933 | 0.0000 | 55.2350 | 92.2853 | 0.3764 |

**Table 9** Equilibrium results of SCN1 associated with markets.

| Market $m$ | $Q_m^E$ | $p_m^E$ |
|---|---|---|
| m1 | 0.0 | 82.9 |
| m2 | 3.2499 | 92.7869 |

*Initialization data of SCN2:*

**Table 10** The parameters of transportation cost function for SCN2.

| No.$(i,j)$ | $a_{ij}$ | $b_{ij}$ | $c_{ij}$ |
|---|---|---|---|
| 1 | 0.0004 | 0.00002 | 0.5 |
| 2 | 0.0004 | 0.00003 | 0.5 |
| 3 | 0.0005 | 0.00002 | 0.5 |
| 4 | 0.0004 | 0.00004 | 0.5 |
| 5 | 0.0004 | 0.00005 | 0.5 |
| 6 | 0.0006 | 0.00002 | 0.5 |
| 7 | 0.0004 | 0.00006 | 0.5 |
| 8 | 0.0008 | 0.00002 | 0.5 |
| 9 | 0.0004 | 0.00006 | 0.5 |
| 10 | 0.0008 | 0.00002 | 0.5 |

*One solution to SCN2 with objective value 0 is given below:*

**Table 11** Equilibrium results of SCN1 associated with links.

| No. $(i,j)$ | $p_i^O$ | $c_{ij}$ | $p_j^E$ | $p_j^E - c_{ij} - p_i^O$ | $f_{i,j}$ |
|---|---|---|---|---|---|
| 1 | 1.0278 | 7.5143 | 8.5421 | 0.0000 | 51.7761 |

| | 2 | 0.3526 | 89.1498 | 89.5024 | 0.0000 | 120.8110 |
|---|---|---|---|---|---|---|
| | 3 | 41.2228 | 25.3640 | 66.5868 | 0.0000 | 74.4062 |
| | 4 | 41.2228 | 18.4301 | 59.6529 | 0.0000 | 70.7426 |
| | 5 | 81.6152 | 0.5826 | 82.1979 | 0.0000 | 11.2797 |
| | 6 | 81.6152 | 10.6798 | 92.2951 | 0.0000 | 63.1265 |
| | 7 | 81.6724 | 0.5254 | 82.1979 | 0.0000 | 7.7615 |
| | 8 | 81.6724 | 10.6226 | 92.2951 | 0.0000 | 62.9812 |
| | 9 | 41.2228 | 40.9751 | 82.1979 | 0.0000 | 99.8638 |
| | 10 | 41.2228 | 40.9750 | 82.1979 | 0.0000 | 100.1616 |

**Table 12** Equilibrium results of SCN1 associated with material suppliers.

| Supplier $s$ | $Q_s^E$ | $Q_s^O$ | $Q_s^H$ | $p_s^O$ | $\lambda_s$ |
|---|---|---|---|---|---|
| s1 | 51.7761 | 51.7761 | 0.0000 | 1.0278 | 0.5130 |
| s2 | 120.8110 | 120.8110 | 0.0000 | 0.3526 | 0.0000 |

**Table 13** Equilibrium results of SCN1 associated with manufacturers.

| Manufacturer $i$ | $Q_i^P$ | $Q_i^O$ | $Q_i^H$ | $Q_{i,s1}^H$ | $Q_{i,s2}^H$ | $p_i^O$ | $p_{i,s1}^E$ | $p_{i,s2}^E$ | $\lambda_i$ |
|---|---|---|---|---|---|---|---|---|---|
| p1 | 345.1743 | 345.1743 | 0.0000 | 0.0000 | 0.0000 | 41.2228 | 8.5421 | 89.5024 | 0.3782 |

**Table 14** Equilibrium results of SCN1 associated with retailers.

| Retailer $j$ | $Q_j^E$ | $Q_j^O$ | $Q_j^H$ | $p_j^E$ | $p_j^O$ | $\lambda_j$ |
|---|---|---|---|---|---|---|
| r1 | 74.4062 | 74.4062 | 0.0000 | 66.5868 | 81.6152 | 0.2231 |
| r2 | 70.7426 | 70.7426 | 0.0000 | 59.6529 | 81.6724 | 0.3626 |

**Table 15** Equilibrium results of SCN1 associated with markets.

| Market $m$ | $Q_m^E$ | $p_m^E$ |
|---|---|---|
| m1 | 219.0666 | 82.1978 |
| m2 | 126.1077 | 92.2950 |

*Initialization data of SCN3:*

**Table 16** Parameters related to costs for SCN2.

| Spots | FC | VC | | HC | | TC | |
|---|---|---|---|---|---|---|---|
| | | $a$ | $b$ | $a$ | $b$ | $a$ | $b$ |
| s1 | 334 | 0.02 | 0.000012 | 0.01 | 0.05 | 0.002 | 0.01 |
| s2 | 411 | 0.01 | 0.000001 | 0.02 | 0.0 | 0.002 | 0.02 |
| s3 | 350 | 0.002 | 0.000002 | 0.01 | 0.0 | 0.001 | 0.02 |
| s4 | 415 | 0.004 | 0.000014 | 0.01 | 0.0 | 0.0002 | 0.01 |
| p1 | 80 | 0.001 | 0.000003 | 0.001 | 0.0 | 0.0001 | 0.0 |
| p2 | 100 | 0.001 | 0.000001 | 0.001 | 0.0 | 0.0002 | 0.0 |
| r1 | 10 | 0.00014 | 0.000001 | 0.005 | 0.0 | 0.005 | 0.0 |
| r2 | 12 | 0.0002 | 0.0000021 | 0.004 | 0.0 | 0.0006 | 0.0 |
| r3 | 0 | 0.0002 | 0.000001 | 0.001 | 0.0 | 0.0001 | 0.0 |

**Table 17** Parameters related to production process and holding cost for SCN2.

| Manufacturer $i$ | $r_{i,1}$ | $r_{i,2}$ | $r_{i,T}$ | HC of material one | | HC of material two | |
|---|---|---|---|---|---|---|---|
| | | | | $a$ | $b$ | $a$ | $b$ |
| p1 | 0.55 | 0.45 | 1.0 | 0.001 | 0.0 | 0.001 | 0.0 |
| p2 | 0.6 | 0.4 | 1.0 | 0.004 | 0.0 | 0.008 | 0.0 |

**Table 18** Parameters associated with pricing function for SCN2.

| Market $m$ | $p_m^{Max}$ | $a_m$ | $b_m$ |
|---|---|---|---|
| m1 | 63..2 | 0.0032 | 0.000076 |
| m2 | 63.8 | 0.004 | 0.000055 |
| m3 | 63.7 | 0.00316 | 0.000068 |
| m4 | 63.5 | 0.0046 | 0.000087 |

**Table 19** The parameters associated with transportation cost function for SCN2.

| No.($i,j$) | $a_{ij}$(E-4) | $b_{ij}$(E-5) | $c_{ij}$ |
|---|---|---|---|
| 1 | 4 | 2 | 0.35 |
| 2 | 4   | 3   | 0.45 |
| 3 | 5   | 2   | 0.35 |
| 4 | 4   | 4   | 0.5  |
| 5 | 4   | 5   | 0.42 |
| 6 | 6   | 2   | 0.35 |
| 7 | 4   | 6   | 0.15 |
| 8 | 4.8 | 4.2 | 0.15 |
| 9 | 5   | 9   | 0.28 |
| 10 | 4  | 4   | 0.45 |
| 11 | 6  | 2   | 0.35 |
| 12 | 4  | 4   | 0.54 |
| 13 | 7  | 2   | 0.25 |
| 14 | 2  | 3   | 0.51 |
| 15 | 4  | 2   | 0.45 |
| 16 | 6  | 3   | 0.25 |
| 17 | 4  | 2   | 0.35 |
| 18 | 4  | 5   | 0.35 |
| 19 | 2  | 2   | 0.45 |
| 20 | 4  | 7   | 0.68 |
| 21 | 5  | 2   | 0.52 |
| 22 | 3  | 8   | 0.53 |
| 23 | 6  | 2   | 0.37 |
| 24 | 4  | 3   | 0.49 |
| 25 | 5  | 3   | 0.45 |
| 26 | 4  | 6   | 0.5  |

*One solution to SCN3 with objective value 0.19268 is given below:*

**Table 20** Equilibrium results of SCN2 associated with links.

| No. $(i,j)$ | $p_i^O$ | $c_{ij}$ | $p_j^E$ | $p_j^E - c_{ij} - p_i^O$ | $f_{i,j}$ |
|---|---|---|---|---|---|
| 1  | 91.1460  | 0.3510 | 91.4861  | -0.0109 | 1.7394  |
| 2  | 91.1460  | 0.4518 | 91.5979  | 0.0000  | 2.4174  |
| 3  | 91.1325  | 0.3536 | 91.4861  | 0.0000  | 3.1273  |
| 4  | 91.1325  | 0.5012 | 91.5979  | -0.0358 | 1.8357  |
| 5  | 120.7559 | 0.4212 | 121.1771 | 0.0000  | 1.9482  |
| 6  | 120.7559 | 0.3510 | 121.1068 | 0.0000  | 1.4016  |
| 7  | 121.0257 | 0.1514 | 121.1771 | 0.0000  | 2.0336  |
| 8  | 121.0257 | 0.1509 | 121.1068 | -0.0697 | 1.4338  |
| 9  | 61.7569  | 0.4579 | 62.2148  | 0.0000  | 13.6961 |
| 10 | 61.7569  | 0.4554 | 62.2123  | 0.0000  | 4.0010  |
| 11 | 61.7569  | 0.3500 | 62.1069  | 0.0000  | 0.0000  |
| 12 | 61.8217  | 0.5400 | 62.2148  | -0.1470 | 0.0540  |
| 13 | 61.8217  | 0.3906 | 62.2123  | 0.0000  | 14.1231 |
| 14 | 61.8217  | 0.5100 | 62.1069  | -0.2248 | 0.0000  |
| 15 | 63.3793  | 0.4500 | 63.2000  | -0.6293 | 0.0000  |
| 16 | 63.3793  | 0.3419 | 63.7212  | 0.0000  | 13.7501 |
| 17 | 63.3793  | 0.3500 | 63.6611  | -0.0682 | 0.0000  |
| 18 | 63.3793  | 0.3500 | 63.5000  | -0.2293 | 0.0000  |
| 19 | 63.0295  | 0.4500 | 63.2000  | -0.2795 | 0.0000  |
| 20 | 63.0295  | 0.6917 | 63.7212  | 0.0000  | 5.8904  |
| 21 | 63.0295  | 0.6316 | 63.6611  | 0.0000  | 12.2337 |
| 22 | 63.0295  | 0.5300 | 63.5000  | -0.0595 | 0.0000  |
| 23 | 65.1579  | 0.3700 | 63.2000  | -2.3279 | 0.0000  |
| 24 | 65.1579  | 0.4900 | 63.7212  | -1.9267 | 0.0000  |
| 25 | 65.1579  | 0.4500 | 63.6611  | -1.9468 | 0.0000  |
| 26 | 65.1579  | 0.5000 | 63.5000  | -2.1579 | 0.0000  |

**Table 21** Equilibrium results of SCN2 associated with material suppliers.

| Supplier $s$ | $Q_s^E$ | $Q_s^O$ | $Q_s^H$ | $p_s^O$ | $\lambda_s$ |
|---|---|---|---|---|---|
| s1 | 4.1567 | 4.1567 | 0.0000 | 91.1460  | 0.1334 |
| s2 | 4.9631 | 4.9631 | 0.0000 | 91.1325  | 0.0990 |
| s3 | 3.3498 | 3.3498 | 0.0000 | 120.7559 | 0.1550 |
| s4 | 3.4674 | 3.4674 | 0.0000 | 121.0257 | 0.0109 |

**Table 22** Equilibrium results of SCN2 associated with manufacturers.

| Manufacturer $i$ | $Q_i^P$ | $Q_i^O$ | $Q_i^H$ | $Q_{i,s1}^H$ | $Q_{i,s2}^H$ | $p_i^O$ | $p_{i,s1}^E$ | $p_{i,s2}^E$ | $\lambda_i$ |
|---|---|---|---|---|---|---|---|---|---|
| p1 | 17.6971 | 17.6971 | 0.0000 | 0.0000 | 0.0000 | 61.7569 | 91.4861 | 121.1771 | 0.0845 |
| p2 | 14.1771 | 14.1771 | 0.0000 | 0.0000 | 0.0000 | 61.8217 | 91.5979 | 121.1068 | 0.0522 |

**Table 23** Equilibrium results of SCN2 associated with retailers.

| Retailer $j$ | $Q_j^E$ | $Q_j^O$ | $Q_j^H$ | $p_j^E$ | $p_j^O$ | $\lambda_j$ |
|---|---|---|---|---|---|---|
| r1 | 1.375E+01 | 1.375E+01 | 0.000E+00 | 6.221E+01 | 6.338E+01 | 6.864E-03 |
| r2 | 1.812E+01 | 1.812E+01 | 0.000E+00 | 6.221E+01 | 6.303E+01 | 2.452E-03 |
| r3 | 4.486E-08 | 4.486E-08 | 0.000E+00 | 6.211E+01 | 6.516E+01 | 4.912E-02 |

**Table 24** Equilibrium results of SCN2 associated with markets.

| Market $m$ | $Q_m^E$ | $p_m^E$ |
|---|---|---|
| m1 | 1.253E-07 | 6.320E+01 |
| m2 | 1.964E+01 | 6.372E+01 |
| m3 | 1.223E+01 | 6.366E+01 |
| m4 | 4.214E-08 | 6.350E+01 |

*Initialization data of SCN4:*

**Table 25** Parameters related to costs for SCN3.

| Spots | FC | VC | | HC | | TC | |
|---|---|---|---|---|---|---|---|
| | | $a$ | $b$ | $a$ | $b$ | $a$ | $b$ |
| s1 | 20 | 0.02 | 0.000012 | 0.0001 | 0.0 | 0.002 | 0.0 |
| s2 | 10 | 0.02 | 0.00001 | 0.0001 | 0.0 | 0.002 | 0.0 |
| s3 | 30 | 0.002 | 0.000002 | 0.0001 | 0.0 | 0.001 | 0.0 |
| s4 | 20 | 0.004 | 0.000014 | 0.0001 | 0.0 | 0.0002 | 0.0 |
| p1 | 50 | 0.003 | 0.000003 | 0.01 | 0.04 | 0.0001 | 0.0 |
| p2 | 80 | 0.001 | 0.000001 | 0.01 | 0.03 | 0.0002 | 0.0 |
| r1 | 10 | 0.0001 | 0.000001 | 0.0005 | 0.0 | 0.005 | 0.0 |
| r2 | 10 | 0.0002 | 0.000001 | 0.04 | 0.0 | 0.0006 | 0.0 |
| r3 | 10 | 0.0002 | 0.000001 | 0.01 | 0.0 | 0.0001 | 0.0 |
| r4 | 20 | 0.0002 | 0.000001 | 0.01 | 0.0 | 0.0001 | 0.0 |

**Table 26** Parameters related to production process and holding cost for SCN3.

| Manufacturer $i$ | $r_{i,1}$ | $r_{i,2}$ | $r_{i,T}$ | HC of material one | | HC of material two | |
|---|---|---|---|---|---|---|---|
| | | | | $a$ | $b$ | $a$ | $b$ |
| p1 | 0.6 | 0.4 | 1.0 | 0.05 | 0.06 | 0.05 | 0.06 |
| p2 | 0.3 | 0.7 | 1.0 | 0.04 | 0.05 | 0.08 | 0.05 |

**Table 27** Parameters associated with pricing function for SCN3.

| Market $m$ | $p_m^{Max}$ | $a_m$ | $b_m$ |
|---|---|---|---|
| m1 | 32.0 | 0.0032 | 0.00076 |
| m2 | 30.8 | 0.005 | 0.00045 |

**Table 28** The parameters associated with transportation cost function for SCN3.

| No.$(i,j)$ | $a_{ij}$(E-4) | $b_{ij}$(E-5) | $c_{ij}$ |
|---|---|---|---|
| 1 | 4 | 3.2 | 0.5 |
| 2 | 4 | 3 | 0.5 |
| 3 | 5 | 2 | 0.5 |
| 4 | 4 | 4 | 0.5 |
| 5 | 4 | 5 | 0.5 |
| 6 | 6 | 2 | 0.5 |
| 7 | 4 | 5 | 0.5 |
| 8 | 8 | 2 | 0.5 |
| 9 | 5 | 9 | 0.5 |
| 10 | 4 | 4 | 0.5 |
| 11 | 6 | 2 | 0.5 |
| 12 | 4 | 4 | 0.5 |
| 13 | 7 | 2 | 0.5 |
| 14 | 2 | 3 | 0.5 |
| 15 | 6.4 | 2 | 0.5 |
| 16 | 6 | 3 | 0.5 |
| 17 | 4 | 2 | 0.5 |
| 18 | 8.4 | 7.5 | 0.5 |
| 19 | 2 | 2 | 0.5 |
| 20 | 4 | 7 | 0.5 |
| 21 | 5 | 2 | 0.5 |
| 22 | 3 | 8 | 0.5 |
| 23 | 6 | 2 | 0.5 |
| 24 | 4 | 5.3 | 0.5 |

*One solution to SCN4 with objective value 0.01190 is given below:*

**Table 29** Equilibrium results of SCN3 associated with links.

| No. $(i,j)$ | $p_i^O$ | $c_{ij}$ | $p_j^E$ | $p_j^E - c_{ij} - p_i^O$ | $f_{i,j}$ |
|---|---|---|---|---|---|
| 1 | 9.9497 | 0.5005 | 10.4502 | 0.0000 | 1.1215 |
| 2 | 9.9497 | 0.5004 | 10.4501 | 0.0000 | 0.9374 |
| 3 | 10.0000 | 0.5000 | 10.4502 | -0.0498 | 0.0000 |
| 4 | 10.0000 | 0.5000 | 10.4501 | -0.0499 | 0.0000 |
| 5 | 10.0000 | 0.5000 | 8.4067 | -2.0933 | 0.0000 |
| 6 | 10.0000 | 0.5000 | 8.4083 | -2.0917 | 0.0000 |
| 7 | 7.9064 | 0.5003 | 8.4067 | 0.0000 | 0.7477 |
| 8 | 7.9064 | 0.5019 | 8.4083 | 0.0000 | 2.1872 |
| 9 | 25.6968 | 0.5000 | 26.1968 | 0.0000 | 0.0892 |
| 10 | 25.6968 | 0.5000 | 26.1961 | -0.0007 | 0.0000 |
| 11 | 25.6968 | 0.5001 | 26.1969 | 0.0000 | 0.1209 |
| 12 | 25.6968 | 0.5021 | 26.1957 | -0.0032 | 3.5283 |
| 13 | 25.6956 | 0.5016 | 26.1968 | -0.0003 | 2.0697 |
| 14 | 25.6956 | 0.5005 | 26.1961 | 0.0000 | 1.8981 |
| 15 | 25.6956 | 0.5014 | 26.1969 | -0.0001 | 2.0096 |
| 16 | 25.6956 | 0.5002 | 26.1957 | 0.0000 | 0.2718 |
| 17 | 31.4646 | 0.5010 | 31.9656 | 0.0000 | 2.1589 |
| 18 | 31.4646 | 0.5000 | 30.8000 | -1.1646 | 0.0000 |
| 19 | 31.4652 | 0.5005 | 31.9656 | 0.0000 | 1.8982 |
| 20 | 31.4652 | 0.5000 | 30.8000 | -1.1652 | 0.0000 |
| 21 | 31.4644 | 0.5012 | 31.9656 | 0.0000 | 2.1306 |
| 22 | 31.4644 | 0.5000 | 30.8000 | -1.1644 | 0.0000 |
| 23 | 31.4628 | 0.5028 | 31.9656 | 0.0000 | 3.8001 |
| 24 | 31.4628 | 0.5000 | 30.8000 | -1.1628 | 0.0000 |

**Table 30** Equilibrium results of SCN3 associated with material suppliers.

| Supplier $s$ | $Q_s^E$ | $Q_s^O$ | $Q_s^H$ | $p_s^O$ | $\lambda_s$ |
|---|---|---|---|---|---|
| s1 | 2.059E+00 | 2.059E+00 | 0.000E+00 | 9.950E+00 | 2.197E-02 |
| s2 | 0.000E+00 | 0.000E+00 | 0.000E+00 | 1.000E+01 | 2.115E-06 |
| s3 | 0.000E+00 | 0.000E+00 | 0.000E+00 | 1.000E+01 | 4.959E-02 |
| s4 | 2.935E+00 | 2.935E+00 | 0.000E+00 | 7.906E+00 | 1.595E-01 |

**Table 31** Equilibrium results of SCN3 associated with manufacturers.

| Manufacturer $i$ | $Q_i^P$ | $Q_i^O$ | $Q_i^H$ | $Q_{i,s1}^H$ | $Q_{i,s2}^H$ | $p_i^O$ | $p_{i,s1}^E$ | $p_{i,s2}^E$ | $\lambda_i$ |
|---|---|---|---|---|---|---|---|---|---|
| p1 | 3.7385 | 3.7385 | 0.0000 | 0.0000 | 0.0000 | 25.6968 | 10.4502 | 8.4067 | 0.4124 |
| p2 | 6.2493 | 6.2493 | 0.0000 | 0.0000 | 0.0000 | 25.6956 | 10.4501 | 8.4083 | 0.4842 |

**Table 32** Equilibrium results of SCN3 associated with retailers.

| Retailer $j$ | $Q_j^E$ | $Q_j^O$ | $Q_j^H$ | $p_j^E$ | $p_j^O$ | $\lambda_j$ |
|---|---|---|---|---|---|---|
| r1 | 2.1589 | 2.1589 | 0.0000 | 26.1968 | 31.4646 | 0.0205 |
| r2 | 1.8982 | 1.8982 | 0.0000 | 26.1961 | 31.4652 | 0.0000 |
| r3 | 2.1306 | 2.1306 | 0.0000 | 26.1969 | 31.4644 | 0.0186 |
| r4 | 3.8001 | 3.8001 | 0.0000 | 26.1957 | 31.4628 | 0.0001 |

**Table 33** Equilibrium results of SCN3 associated with markets.

| Market $m$ | $Q_m^E$ | $p_m^E$ |
|---|---|---|
| m1 | 9.9877 | 31.9656 |
| m2 | 6.953E-07 | 30.8000 |

*Initialization data of SCN5:*

**Table 34** Parameters related to costs for SCN3.

| Spots | FC | VC | | HC | | TC | |
|---|---|---|---|---|---|---|---|
| | | $a$ | $b$ | $a$ | $b$ | $a$ | $b$ |
| s1 | 10 | 0.02 | 0.000012 | 0.0001 | 0.0 | 0.002 | 0.0 |
| s2 | 15 | 0.02 | 0.00001 | 0.0001 | 0.0 | 0.002 | 0.0 |
| s3 | 20 | 0.002 | 0.000002 | 0.0001 | 0.0 | 0.001 | 0.0 |
| s4 | 22 | 0.004 | 0.000014 | 0.0001 | 0.0 | 0.0002 | 0.0 |
| p1 | 50 | 0.003 | 0.000003 | 0.01 | 0.04 | 0.0001 | 0.0 |
| p2 | 80 | 0.001 | 0.000001 | 0.01 | 0.03 | 0.0002 | 0.0 |
| w1 | 1 | 0.00005 | 0.0000001 | 0.00001 | 0.0 | 0.002 | 0.0 |
| r1 | 10 | 0.0001 | 0.000001 | 0.0005 | 0.0 | 0.005 | 0.0 |
| r2 | 10 | 0.0002 | 0.000001 | 0.04 | 0.0 | 0.0006 | 0.0 |
| r3 | 10 | 0.0002 | 0.000001 | 0.01 | 0.0 | 0.0001 | 0.0 |
| r4 | 20 | 0.0002 | 0.000001 | 0.01 | 0.0 | 0.0001 | 0.0 |

Table 35 Parameters related to production process and holding cost for SCN3.

| Manufacturer $i$ | $r_{i,1}$ | $r_{i,2}$ | $r_{i,T}$ | HC of material one | | HC of material two | |
|---|---|---|---|---|---|---|---|
| | | | | a | b | a | b |
| p1 | 0.6 | 0.4 | 1.0 | 0.05 | 0.06 | 0.05 | 0.06 |
| p2 | 0.3 | 0.7 | 1.0 | 0.04 | 0.05 | 0.08 | 0.05 |

Table 36 Parameters associated with pricing function for SCN3.

| Market $m$ | $p_m^{Max}$ | $a_m$ | $b_m$ |
|---|---|---|---|
| m1 | 33.0 | 0.0032 | 0.00066 |
| m2 | 32.8 | 0.005 | 0.00055 |

Table 37 The parameters associated with transportation cost function for SCN3.

| No.$(i,j)$ | $a_{ij}$(E-4) | $b_{ij}$(E-5) | $c_{ij}$ |
|---|---|---|---|
| 1 | 4 | 3.2 | 0.5 |
| 2 | 4 | 3 | 0.5 |
| 3 | 5 | 2 | 0.5 |
| 4 | 4 | 4 | 0.5 |
| 5 | 4 | 5 | 0.5 |
| 6 | 6 | 2 | 0.5 |
| 7 | 4 | 5 | 0.5 |
| 8 | 8 | 2 | 0.5 |
| 13 | 7.7 | 2 | 0.5 |
| 14 | 5.2 | 3 | 0.5 |
| 15 | 8.4 | 2 | 0.5 |
| 16 | 7.6 | 3 | 0.5 |
| 17 | 4 | 2 | 0.5 |
| 18 | 8.4 | 7.5 | 0.5 |
| 19 | 2 | 2 | 0.5 |
| 20 | 4 | 7 | 0.5 |
| 21 | 5 | 2 | 0.5 |
| 22 | 3 | 8 | 0.5 |
| 23 | 6 | 2 | 0.5 |
| 24 | 4 | 5.3 | 0.5 |
| 25 | 0.4 | 0.5 | 0.15 |
| 26 | 0.6 | 0.2 | 0.15 |
| 27 | 0.4 | 0.5 | 0.15 |
| 28 | 0.8 | 0.2 | 0.15 |
| 29 | 0.5 | 0.9 | 0.15 |
| 30 | 0.4 | 0.4 | 0.51 |

*One solution to SCN5 with objective value 0.14627 is given below:*

Table 38 Equilibrium results of SCN3 associated with links.

| No.$(i,j)$ | $p_i^O$ | $C_{ij}$ | $p_j^E$ | $p_j^E - C_{ij} - p_i^O$ | $f_{i,j}$ |
|---|---|---|---|---|---|
| 1 | 10.0000 | 0.5000 | 8.1206 | -2.3794 | 0.0000 |
| 2 | 10.0000 | 0.5000 | 8.1196 | -2.3804 | 0.0000 |
| 3 | 7.6190 | 0.5016 | 8.1206 | 0.0000 | 2.6460 |
| 4 | 7.6190 | 0.5006 | 8.1196 | 0.0000 | 1.3133 |
| 5 | 12.0617 | 0.5007 | 12.5614 | -0.0011 | 1.4491 |
| 6 | 12.0617 | 0.5009 | 12.5626 | 0.0000 | 1.3485 |
| 7 | 12.0612 | 0.5001 | 12.5614 | 0.0000 | 0.3150 |
| 8 | 12.0612 | 0.5015 | 12.5626 | -0.0001 | 1.7159 |
| 13 | 18.0749 | 0.5000 | 18.5749 | 0.0000 | 0.0420 |
| 14 | 18.0749 | 0.5039 | 18.5776 | -0.0012 | 4.9353 |
| 15 | 18.0749 | 0.5025 | 18.5774 | 0.0000 | 2.6656 |
| 16 | 18.0749 | 0.5009 | 18.5758 | 0.0000 | 1.1122 |
| 17 | 32.4397 | 0.5014 | 32.9410 | 0.0000 | 2.7866 |
| 18 | 32.4397 | 0.5000 | 32.7999 | -0.1398 | 0.0001 |
| 19 | 32.4383 | 0.5029 | 32.9410 | -0.0002 | 5.8743 |
| 20 | 32.4383 | 0.5000 | 32.7999 | -0.1385 | 0.0125 |
| 21 | 32.4342 | 0.5069 | 32.9410 | -0.0001 | 7.4552 |
| 22 | 32.4342 | 0.5000 | 32.7999 | -0.1343 | 0.0000 |
| 23 | 32.4402 | 0.5009 | 32.9410 | -0.0001 | 1.4462 |
| 24 | 32.4402 | 0.5000 | 32.7999 | -0.1403 | 0.0008 |
| 25 | 13.9672 | 0.1542 | 14.1214 | 0.0000 | 8.8202 |
| 26 | 18.0749 | 0.1500 | 14.1214 | -4.1036 | 0.0004 |
| 27 | 18.4275 | 0.1503 | 18.5749 | -0.0028 | 2.7446 |
| 28 | 18.4275 | 0.1501 | 18.5776 | 0.0000 | 0.9515 |
| 29 | 18.4275 | 0.1510 | 18.5774 | -0.0011 | 4.7896 |
| 30 | 18.4275 | 0.5100 | 18.5758 | -0.3617 | 0.3348 |

**Table 39** Equilibrium results of SCN3 associated with material suppliers.

| Supplier $s$ | $Q_s^E$ | $Q_s^O$ | $Q_s^H$ | $p_s^O$ | $\lambda_s$ |
|---|---|---|---|---|---|
| s1 | 0.0000 | 0.0000 | 0.0000 | 10.0000 | 0.0578 |
| s2 | 3.9594 | 3.9594 | 0.0000 | 7.6190 | 0.9995 |
| s3 | 2.7976 | 2.7976 | 0.0000 | 12.0617 | 0.6865 |
| s4 | 2.0309 | 2.0309 | 0.0000 | 12.0612 | 0.1130 |

**Table 40** Equilibrium results of SCN3 associated with manufacturers.

| Manufacturer $i$ | $Q_i^P$ | $Q_i^O$ | $Q_i^H$ | $Q_{i,s1}^H$ | $Q_{i,s2}^H$ | $p_i^O$ | $p_{i,s1}^E$ | $p_{i,s2}^E$ | $\lambda_i$ |
|---|---|---|---|---|---|---|---|---|---|
| p1 | 8.8202 | 8.8202 | 0.0000 | 0.0000 | 0.0000 | 13.9672 | 8.1206 | 12.5614 | 0.3151 |
| p2 | 8.7555 | 8.7555 | 0.0000 | 0.0000 | 0.0000 | 18.0749 | 8.1196 | 12.5626 | 0.2252 |

**Table 41** Equilibrium results of SCN3 associated with wholesaler and retailers.

| $j$ | $Q_j^E$ | $Q_j^O$ | $Q_j^H$ | $p_j^E$ | $p_j^O$ | $\lambda_j$ |
|---|---|---|---|---|---|---|
| w1 | 8.8206 | 8.8206 | 0.0000 | 14.1214 | 18.4275 | 0.2944 |
| r1 | 2.7867 | 2.7867 | 0.0000 | 18.5749 | 32.4397 | 0.4633 |
| r2 | 5.8868 | 5.8868 | 0.0000 | 18.5776 | 32.4383 | 0.5998 |
| r3 | 7.4552 | 7.4552 | 0.0000 | 18.5774 | 32.4342 | 0.6283 |
| r4 | 1.4470 | 1.4470 | 0.0000 | 18.5758 | 32.4402 | 0.0013 |

**Table 42** Equilibrium results of SCN3 associated with markets.

| Market $m$ | $Q_m^E$ | $p_m^E$ |
|---|---|---|
| m1 | 17.5623 | 32.9410 |
| m2 | 0.0135 | 32.7999 |

References:

Afshar, A., Haddad, O.B., Marino, M.A., Adams, B.J., 2007. Honey-bee mating optimization (HBMO) algorithm for optimal reservoir operation. Journal of the Franklin Institute 344, 452-462.

Atashpaz-Gargari, E., Lucas, C., 2007. Imperialist competitive algorithm: an algorithm for optimization inspired by imperialistic competition, 2007 IEEE congress on evolutionary computation. Ieee, pp. 4661-4667.

Auger, A., Hansen, N., 2005. A restart CMA evolution strategy with increasing population size, 2005 IEEE Congress on Evolutionary Computation, pp. 1769-1776 Vol. 1762.

Beyer, H.-G., Sendhoff, B., 2008. Covariance matrix adaptation revisited–the CMSA evolution strategy–, International Conference on Parallel Problem Solving from Nature. Springer, pp. 123-132.

Brandao, M.S., Godinho-Filho, M., 2022. Is a multiple supply chain management perspective a new way to manage global supply chains toward sustainability? Journal of Cleaner Production 375, 134046.

Carvalho, H., Naghshineh, B., Govindan, K., Cruz-Machado, V., 2022. The resilience of on-time delivery to capacity and material shortages: An empirical investigation in the automotive supply chain. Computers & Industrial Engineering 171, 108375.

Chan, C.K., Zhou, Y., Wong, K.H., 2019. An equilibrium model of the supply chain network under multi-attribute behaviors analysis. European Journal of Operational Research 275, 514-535.

Chen, J., Xu, Z., Huang, D., Fang, C., Wang, X., Zhang, J., 2020. Automotive supply chain networks equilibrium model under uncertain payment delay and objective weights. Computers & Industrial Engineering 150, 106866.

Cheng, M.-Y., Prayogo, D., 2014. Symbiotic Organisms Search: A new metaheuristic optimization algorithm. Computers & Structures 139, 98-112.

Civicioglu, P., 2013. Backtracking search optimization algorithm for numerical optimization problems. Applied Mathematics and computation 219, 8121-8144.

Daultani, Y., Kumar, S., Vaidya, O.S., Tiwari, M.K., 2015. A supply chain network equilibrium model for operational and opportunism risk mitigation. International Journal of Production Research 53, 5685-5715.

Dong, J., Zhang, D., Nagurney, A., 2004. A supply chain network equilibrium model with random demands. European Journal of Operational Research 156, 194-212.

Dong, J., Zhang, D., Yan, H., Nagurney, A., 2005. Multitiered Supply Chain Networks: Multicriteria Decision—Making Under Uncertainty. Annals of Operations Research 135, 155-178.

Dorigo, M., Birattari, M., Stutzle, T., 2006. Ant colony optimization. IEEE computational intelligence magazine 1, 28-39.

Erol, O.K., Eksin, I., 2006. A new optimization method: Big Bang–Big Crunch. Advances in Engineering Software 37, 106-111.

Eskandar, H., Sadollah, A., Bahreininejad, A., Hamdi, M., 2012. Water cycle algorithm – A novel metaheuristic optimization method for solving constrained engineering optimization problems. Computers & Structures 110-111, 151-166.

Eusuff, M., Lansey, K., Pasha, F., 2006. Shuffled frog-leaping algorithm: a memetic meta-heuristic for discrete optimization. Engineering optimization 38, 129-154.

Faramarzi, A., Heidarinejad, M., Stephens, B., Mirjalili, S., 2020. Equilibrium optimizer: A novel optimization algorithm. Knowledge-based systems 191, 105190.

Feng, P., Zhou, X., Zhang, D., Chen, Z., Wang, S., 2022. The impact of trade policy on global supply chain network equilibrium: A new perspective of product-market chain competition. Omega 109, 102612.

Feng, Z.-F., Wang, Z.-p., Chen, Y., 2014. The equilibrium of closed-loop supply chain supernetwork with time-dependent parameters. Transportation Research Part E-logistics and Transportation Review 64, 1-11.

Formato, R., 2007. Central force optimization: a new metaheuristic with applications in applied electromagnetics. Progress In Electromagnetics Research 77, 425-491.

Gandomi, A.H., Yang, X.-S., Alavi, A.H., Talatahari, S., 2013. Bat algorithm for constrained optimization tasks. Neural Computing and Applications 22, 1239-1255.

Glover, F., 1990. Tabu search: A tutorial. Interfaces 20, 74-94.

Gupta, R., Goswami, M., Daultani, Y., Biswas, B., Allada, V., 2023. Profitability and pricing decision-making structures in presence of uncertain demand and green technology investment for a three tier supply chain. Computers & Industrial Engineering 179, 109190.

Hammond, D., Beullens, P., 2007. Closed-loop supply chain network equilibrium under legislation. European Journal of Operational Research 183, 895-908.

Hansen, N., Ostermeier, A., 1996. Adapting arbitrary normal mutation distributions in evolution strategies: The covariance matrix adaptation, Proceedings of IEEE international conference on evolutionary computation. IEEE, pp. 312-317.

Hashim, F.A., Hussien, A.G., 2022. Snake Optimizer: A novel meta-heuristic optimization algorithm. Knowledge-Based Systems 242, 108320.

Hatamlou, A., 2013. Black hole: A new heuristic optimization approach for data clustering. Information sciences 222, 175-184.

He, P., Zhang, G., Wang, T.-Y., Si, Y., 2023. Optimal two-period pricing strategies in a dual-channel supply chain considering market change. Computers & Industrial Engineering 179, 109193.

He, S.-X., 2023. Truss optimization with frequency constraints using the medalist learning algorithm. Structures 55, 1-15.

He, S.-X., Cui, Y.-T., 2023a. Medalist learning algorithm for configuration optimization of trusses. Applied Soft Computing 148, 110889.

He, S.-X., Cui, Y.-T., 2023b. Multiscale medalist learning algorithm and its application in engineering. Acta Mechanica.

He, S.-X., Cui, Y.-T., 2023c. A novel variational inequality approach for modeling the optimal equilibrium in multi-tiered supply chain networks. Supply Chain Analytics 4, 100039.

Hirano, T., Narushima, Y., 2019. Robust Supply Chain Network Equilibrium Model. Transportation Science 53, 1196-1212.

Holland, J.H., 1992. Genetic algorithms. Scientific american 267, 66-73.

Hsueh, C.-F., Chang, M.-S., 2008. Equilibrium analysis and corporate social responsibility for supply chain integration. European Journal of Operational Research 190, 116-129.

Huggins, E.L., Olsen, T.L., 2003. Supply Chain Management with Guaranteed Delivery. Management Science 49, 1154-1167.

Iacca, G., dos Santos Junior, V.C., de Melo, V.V., 2021. An improved Jaya optimization algorithm with Lévy flight. Expert Systems with Applications 165, 113902.

Jabbarzadeh, A., Haughton, M., Khosrojerdi, A., 2018. Closed-loop supply chain network design under disruption risks: A robust approach with real world application. Computers & Industrial Engineering 116, 178-191.

Jingqiao, Z., Sanderson, A.C., 2007. JADE: Self-adaptive differential evolution with fast and reliable convergence performance, 2007 IEEE Congress on Evolutionary Computation, pp. 2251-2258.

Karaboga, D., 2005. An idea based on honey bee swarm for numerical optimization. Technical report-tr06, Erciyes university, engineering faculty, computer ….

Kashan, A.H., 2009. League championship algorithm: a new algorithm for numerical function optimization, 2009 international conference of soft computing and pattern recognition. IEEE, pp. 43-48.

Kaveh, A., Bakhshpoori, T., 2016. Water Evaporation Optimization: A novel physically inspired optimization algorithm. Computers & Structures 167, 69-85.

Kaveh, A., Dadras, A., 2017. A novel meta-heuristic optimization algorithm: Thermal exchange optimization. Advances in Engineering Software 110, 69-84.

Kaveh, A., Farhoudi, N., 2013. A new optimization method: Dolphin echolocation. Advances in Engineering Software 59, 53-70.

Kaveh, A., Ghazaan, M.I., 2014. Enhanced colliding bodies optimization for design problems with continuous and discrete variables. Advances in Engineering Software 77, 66-75.

Kaveh, A., Ilchi Ghazaan, M., 2017. A new meta-heuristic algorithm: vibrating particles system. Scientia Iranica 24, 551-566.

Kaveh, A., Khayatazad, M., 2012. A new meta-heuristic method: ray optimization. Computers & structures 112, 283-294.

Kaveh, A., Talatahari, S., 2010. A novel heuristic optimization method: charged system search. Acta Mechanica 213, 267-289.

Kaveh, A., Zolghadr, A., 2016. A novel meta-heuristic algorithm: tug of war optimization.

Kennedy, J., Eberhart, R., 1995. Particle swarm optimization, Proceedings of ICNN'95 - International Conference on Neural Networks, pp. 1942-1948 vol.1944.

Khatri, A., Gaba, A., Rana, K., Kumar, V., 2020. A novel life choice-based optimizer. Soft computing 24, 9121-9141.

Kirkpatrick, S., D, G.C., P, V.M., 1983. Simulated annealing. Science 220, 671-680.

Korpeoglu, C.G., Körpeoğlu, E., Cho, S.-H., 2020. Supply Chain Competition: A Market Game Approach. Management Science 66, 5648-5664.

Lee, K.S., Geem, Z.W., 2005. A new meta-heuristic algorithm for continuous engineering optimization: harmony search theory and practice. Computer Methods in Applied Mechanics and Engineering 194, 3902-3933.

Liu, Z., Cruz, J.M., 2012. Supply chain networks with corporate financial risks and trade credits under economic uncertainty. International Journal of Production Economics 137, 55-67.

Liu, Z., Nagurney, A., 2011. Supply Chain Outsourcing Under Exchange Rate Risk and Competition. Omega 39, 539-549.

Liu, Z., Nagurney, A., 2013. Supply chain networks with global outsourcing and quick-response production under demand and cost uncertainty. Annals of Operations Research 208, 251-289.

Liu, Z., Wang, J., 2019. Supply chain network equilibrium with strategic financial hedging using futures. European Journal of Operational Research 272, 962-978.

Ma, J., Zhang, D., Dong, J., Tu, Y., 2020. A supply chain network economic model with time-based competition. European Journal of Operational Research 280, 889-908.

Masoumi, A.H., Yu, M., Nagurney, A., 2017. Mergers and acquisitions in blood banking systems: A supply chain network approach. International Journal of Production Economics 193, 406-421.

Mehrabian, A.R., Lucas, C., 2006. A novel numerical optimization algorithm inspired from weed colonization. Ecological informatics 1, 355-366.

Meng, Z., Pan, J.-S., 2016. Monkey king evolution: a new memetic evolutionary algorithm and its application in vehicle fuel consumption optimization. Knowledge-Based Systems 97, 144-157.

Mirjalili, S., 2015a. The ant lion optimizer. Advances in engineering software 83, 80-98.

Mirjalili, S., 2015b. Moth-flame optimization algorithm: A novel nature-inspired heuristic paradigm. Knowledge-Based Systems 89, 228-249.

Mirjalili, S., 2016a. Dragonfly algorithm: a new meta-heuristic optimization technique for solving single-objective, discrete, and multi-objective problems. Neural computing and applications 27, 1053-1073.

Mirjalili, S., 2016b. SCA: a sine cosine algorithm for solving optimization problems. Knowledge-based systems 96, 120-133.

Mirjalili, S., Gandomi, A.H., Mirjalili, S.Z., Saremi, S., Faris, H., Mirjalili, S.M., 2017. Salp Swarm Algorithm: A bio-inspired optimizer for engineering design problems. Advances in Engineering Software 114, 163-191.

Mirjalili, S., Lewis, A., 2016. The Whale Optimization Algorithm. Advances in Engineering Software 95, 51-67.

Mirjalili, S., Mirjalili, S.M., Hatamlou, A., 2016. Multi-Verse Optimizer: a nature-inspired algorithm for global optimization. Neural Computing and Applications 27, 495-513.

Mirjalili, S., Mirjalili, S.M., Lewis, A., 2014. Grey Wolf Optimizer. Advances in Engineering Software 69, 46-61.

Mirjalili, S.Z., Mirjalili, S., Saremi, S., Faris, H., Aljarah, I., 2018. Grasshopper optimization algorithm for multi-objective optimization problems. Applied Intelligence 48, 805-820.

Mohamed, A.W., Hadi, A.A., Mohamed, A.K., 2020. Gaining-sharing knowledge based algorithm for solving optimization problems: a novel nature-inspired algorithm. International Journal of Machine Learning and Cybernetics 11, 1501-1529.

Nagurney, A., 2021. Supply chain game theory network modeling under labor constraints: Applications to the Covid-19 pandemic. European Journal of Operational Research 293, 880-891.

Nagurney, A., Besik, D., Yu, M., 2018. Dynamics of quality as a strategic variable in complex food supply chain network competition: The case of fresh produce. Chaos: An Interdisciplinary Journal of Nonlinear Science 28, 043124.

Nagurney, A., Cruz, J., Dong, J., Zhang, D., 2005. Supply chain networks, electronic commerce, and supply side and demand side risk. European Journal of Operational Research 164, 120-142.

Nagurney, A., Daniele, P., 2021. International human migration networks under regulations. European Journal of Operational Research 291, 894-905.

Nagurney, A., Ke, K., 2006. Financial networks with intermediation: Risk management with variable weights. European Journal of Operational Research 172, 40-63.

Nagurney, A., Liu, Z., Cojocaru, M.-G., Daniele, P., 2007. Dynamic electric power supply chains and transportation networks: An evolutionary variational inequality formulation. Transportation Research Part E: Logistics and Transportation Review 43, 624-646.

Nagurney, A., Loo, J., Dong, J., Zhang, D., 2002. Supply Chain Networks and Electronic Commerce: A Theoretical Perspective. Netnomics 4, 187-220.

Nagurney, A., Yu, M., 2012. Sustainable fashion supply chain management under oligopolistic competition and brand differentiation. International Journal of Production Economics 135, 532-540.

Nguyen, T., Hoang, B., Nguyen, G., Nguyen, B.M., 2020. A new workload prediction model using extreme learning machine and enhanced tug of war optimization. Procedia Computer Science 170, 362-369.

Oladejo, S.O., Ekwe, S.O., Mirjalili, S., 2024. The Hiking Optimization Algorithm: A novel human-based metaheuristic approach. Knowledge-Based Systems 296, 111880.

Rahmani, R., Yusof, R., 2014. A new simple, fast and efficient algorithm for global optimization over continuous search-space problems: radial movement optimization. Applied Mathematics and Computation 248, 287-300.

Rao, R.V., 2019. Jaya: an advanced optimization algorithm and its engineering applications.

Rao, R.V., Savsani, V.J., Vakharia, D.P., 2011. Teaching–learning-based optimization: a novel method for constrained mechanical design optimization problems. Computer-aided design 43, 303-315.

Rashedi, E., Nezamabadi-pour, H., Saryazdi, S., 2009. GSA: A Gravitational Search Algorithm. Information Sciences 179, 2232-2248.

Sadeeq, H.T., Abdulazeez, A.M., 2022. Giant trevally optimizer (GTO): A novel metaheuristic algorithm for global optimization and challenging engineering problems. Ieee Access 10, 121615-121640.

Sadollah, A., Bahreininejad, A., Eskandar, H., Hamdi, M., 2013. Mine blast algorithm: A new population based algorithm for solving constrained engineering optimization problems. Applied Soft Computing 13, 2592-2612.

Salcedo-Sanz, S., Del Ser, J., Landa-Torres, I., Gil-López, S., Portilla-Figueras, J., 2014. The coral reefs optimization algorithm: a novel metaheuristic for efficiently solving optimization problems. The Scientific World Journal 2014, 739768.

Salhi, A., Fraga, E.S., 2011. Nature-inspired optimisation approaches and the new plant propagation algorithm.

Salimi, H., 2015. Stochastic Fractal Search: A powerful metaheuristic algorithm. Knowledge-Based Systems 75, 1-18.

Sawik, T., 2023. Reshore or not Reshore: A Stochastic Programming Approach to Supply Chain Optimization. Omega 118, 102863.

Simon, D., 2008. Biogeography-Based Optimization. IEEE Transactions on Evolutionary Computation 12, 702-713.

Storn, R., Price, K., 1997. Differential Evolution – A Simple and Efficient Heuristic for global Optimization over Continuous Spaces. Journal of Global Optimization 11, 341-359.

Tanabe, R., Fukunaga, A., 2013. Success-history based parameter adaptation for differential evolution, 2013 IEEE congress on evolutionary computation. IEEE, pp. 71-78.

Tanabe, R., Fukunaga, A.S., 2014. Improving the search performance of SHADE using linear population size reduction, 2014 IEEE Congress on Evolutionary Computation (CEC), pp. 1658-1665.

Van Thieu, N., Mirjalili, S., 2023. MEALPY: An open-source library for latest meta-heuristic algorithms in Python. Journal of Systems Architecture 139, 102871.

Wang, J., Zhou, H., Sun, X., Yuan, Y., 2023. A novel supply chain network evolving model under random and targeted disruptions. Chaos, Solitons & Fractals 170, 113371.

Wu, K., Nagurney, A., Liu, Z., Stranlund, J.K., 2006. Modeling generator power plant portfolios and pollution taxes in electric power supply chain networks: A transportation network equilibrium transformation. Transportation Research Part D: Transport and Environment 11, 171-190.

Xiao, Y.-X., Zhang, R.-Q., 2023. Supply chain network equilibrium considering coordination between after-sale service and product quality. Computers & Industrial Engineering 175, 108848.

Yadav, A., Sadollah, A., Yadav, N., Kim, J.H., 2020. Self-adaptive global mine blast algorithm for numerical optimization. Neural Computing and Applications 32, 2423-2444.

Yang, X.-S., 2010a. Firefly algorithm, stochastic test functions and design optimisation, arXiv preprint arXiv:1003.1409, pp. 1-12.

Yang, X.-S., 2010b. A new metaheuristic bat-inspired algorithm, Nature inspired cooperative strategies for optimization (NICSO 2010). Springer, pp. 65-74.

Yang, X.-S., 2012. Flower pollination algorithm for global optimization, International conference on unconventional computing and natural computation. Springer, pp. 240-249.

Yang, X.-S., Deb, S., 2009. Cuckoo search via Lévy flights, 2009 World congress on nature & biologically inspired computing (NaBIC). Ieee, pp. 210-214.

Yu, M., Nagurney, A., 2013. Competitive food supply chain networks with application to fresh produce. European Journal of Operational Research 224, 273-282.

Zhang, D., 2006. A network economic model for supply chain versus supply chain competition. Omega 34, 283-295.

Zhang, G., Dai, G., Sun, H., Zhang, G., Yang, Z., 2020. Equilibrium in supply chain network with competition and service level between channels considering consumers' channel preferences. Journal of Retailing and Consumer Services 57, 102199.

Zhang, L., Zhou, Y., 2012. A new approach to supply chain network equilibrium models. Computers & Industrial Engineering 63, 82-88.

Zhou, Y., Chan, C.K., Wong, K.H., 2018. A multi-period supply chain network equilibrium model considering retailers' uncertain demands and dynamic loss-averse behaviors. Transportation Research Part E: Logistics and Transportation Review 118, 51-76.